\documentclass{amsart}
\usepackage[utf8]{inputenc}
\usepackage{amsmath}
 \usepackage{stmaryrd}
\usepackage{amssymb,latexsym}
\usepackage{amsfonts,mathrsfs}
\usepackage[margin=1.8in]{geometry}
\usepackage[all,cmtip]{xy}
\usepackage[colorlinks= true,pdfborder={0 0 0}]{hyperref}
\usepackage{graphicx}
 \usepackage{esint}
 \usepackage {verbatim}
\usepackage{color}
\usepackage{enumitem}
\def\pd#1{\dfrac{\partial}{\partial #1}}
\def\f#1#2{\frac{#1}{#2}}

\def\pa{\partial}

\def\n{\nabla}

\def\({\left (}
\def\){\right )}
\def\<{\langle}
\def\>{\rangle}

\newcommand{\bel}[1]{\begin{equation}\label{#1}}

  \newcommand{\beq}{\begin{equation}}

\newcommand{\ba}{\begin{eqnarray}}
\newcommand{\ea}{\end{eqnarray}}

\newcommand{\qe}{\end{equation}}
\newcommand{\eeq}{\end{equation}}

\newtheorem{thm}{Theorem}[section]

\newtheorem{lem}[thm]{Lemma}
\newtheorem{prop}[thm]{Proposition}
\newtheorem{defn}[thm]{Definition}

\newtheorem*{acknowledge}{Acknowledgement}

\newcommand{\norm}[1]{\left\Vert#1\right\Vert}
\newcommand{\abs}[1]{\left\vert#1\right\vert}

\newcommand{\R}{\mathbb R}

\newcommand{\red}[1]{ {\color{red} #1} }
\newcommand{\blue}[1]{ {\color{blue} #1} }

\title[uniqueness of HCMC]{Constant harmonic mean curvature  foliation in asymptotically Schwarzschild spaces-II}
\author{Yaoting Gui, Yuqiao Li, Jun Sun}
\address{Yaoting Gui, School of Mathematical Sciences, Xiamen University, Xiamen, 361005, P. R. China}
\address{Yuqiao Li, Department of Mathematics, Hefei University of Technology, Hefei, 230009, P. R. China}
\address{Jun Sun, School of Mathematics and Statistics, Wuhan University, Wuhan, 430072, P. R. China}
\email{ytgui@bicmr.pku.edu.cn, lyq112@mail.ustc.edu.cn, sunjun@whu.edu.cn}
\begin{document}

\begin{abstract}
    This paper extends the results of \cite{Guilisun}, where the existence of a constant harmonic mean curvature foliation was established in the setting of a 3-dimensional asymptotically Schwarzschild manifold. Here, we generalize this construction to higher dimensions, proving the existence of foliations by constant harmonic mean curvature hypersurfaces in an asymptotically Schwarzschild manifold of arbitrary dimension. Furthermore, in 3 dimensional case, we demonstrate the local uniqueness of this foliation under a stronger decay conditions on the asymptotically Schwarzschild metric.
\end{abstract}

\maketitle

\section{Introduction}

\allowdisplaybreaks
\vspace{.1in}

The existence of foliation by special surfaces in asymptotically Schwarzschild spaces constitutes a fundamental problem in general relativity. 
In their seminal paper, Huisken-Yau (\cite{Huisken1996DefinitionOC}) established the existence of a foliation by stable constant mean curvature (CMC) surfaces in an asymptotically Schwarzschild space with $m>0$. 
These surfaces arise from the volume preserving mean curvature flow introduced by Huisken (\cite{Huisken1987TheVP}).  Moreover, Huisken and Yau proved that the foliation is uniquely determined outside a sufficiently large compact set—dependent on the mean curvature—and used it to define a geometric center of mass.
The uniqueness was later strengthened by Qing-Tian (\cite{QingTian}), who showed that the foliation is, in fact, unique outside a fixed compact set, a property now referred to as global uniqueness. 
Furthermore, Corvino-Wu demonstrated that the geometric center of mass associated with this CMC foliation coincides with the ADM center of mass, provided that the metric is conformally flat at infinity \cite{Corvino}. 
Huang subsequently removed this conformal flatness assumption, thereby generalizing the result to broader settings in \cite{Huang09}.

Recently, we proved the existence of a foliation by constant harmonic mean curvature (CHMC) surfaces in asymptotically Schwarzschild 3-manifolds in (\cite{Guilisun}). 
Our approach utilizes the volume-preserving harmonic mean curvature flow to generate such surfaces and demonstrates that they naturally form a foliation. As a direct consequence, we introduced a new notion of center of mass $C_{HM}$
  and established its equivalence with the ADM mass. However, the extension of these results to higher dimensions, as well as the uniqueness problem remains open and deserving further investigation.

This paper focuses on asymptotically Schwarzschild spaces of general dimension $n\geq3$. 
Specifically, we consider a complete Riemannian $n$-manifold $(N, \bar{g})$  diffeomorphic to $\mathbb{R}^n\backslash B_1(0)$ equipped with a metric $\bar{g}$ that is asymptotically Schwarzschild.
More precisely, the metric takes the form
\begin{equation}\label{metr}
\bar{g}_{\alpha\beta}=\left(1+\frac{m}{2r^{n-2}}\right)^{\frac{4}{n-2}}\delta_{\alpha\beta}+P_{\alpha\beta}, 
\end{equation} 
where $r$ denotes the Euclidean distance, $m>0$  represents the ADM mass and the perturbation term $P_{\alpha\beta}$ satisfies
\[ |\partial^lP_{\alpha\beta}|\leq C_{l+1}r^{1-n-l-\delta}, \quad 0\leq l\leq5, \]
with $\partial$ denoting the partial derivatives with respect to the Euclidean metric $\delta_{\alpha\beta}$ and $\delta\geq 0$ being a constant. 
For notational simplicity, we in the sequel set $c_0=\max(1, m, C_1, C_2, C_3, C_4, C_5, C_6)$.
Note that when $\delta>0$, the decay rate imposed here is stronger than that required in the work of Huisken-Yau in \cite{Huisken1996DefinitionOC}.

Our main theorem establishes the existence of a foliation by constant harmonic mean curvature (CHMC) hypersurfaces for all $\delta\geq 0$. This foliation enables us to define a geometric center of mass, which we prove to be equivalent to the ADM center of mass. Furthermore, under the stronger decay condition $\delta>0$ and in dimension $n=3$, we establish the local uniqueness of this foliation.

\begin{thm}
    Let $(N, \bar{g})$ be an asymptotically Schwarzschild $n$-manifold with metric \eqref{metr}. There is $\sigma_0$ depending only on $c_0$ such that for all $\sigma\geq\sigma_0$, the constant harmonic mean curvature hypersurfaces $\Sigma_{\sigma}$ constructed in Theorem \ref{main} constitute a proper foliation of $N\backslash B_{\sigma_0}(0)$.
    Moreover, the geometric center of mass defined by the CHMC foliation in Definition \ref{center} is equivalent to the ADM center of mass.
    In particular, for $n=3$ and $\delta>0$, $\Sigma_{\sigma}$ constructed in Theorem \ref{main} is the only surface with constant harmonic mean curvature $F_{\sigma}$ contained in $\mathcal{B}_{\sigma} (B_1, B_2, B_3)$, where $\mathcal{B}_{\sigma} (B_1, B_2, B_3)$ is defined in Section 2.
\end{thm}

Following the approach in \cite{Guilisun}, we establish the existence of the constant harmonic mean curvature (CHMC) surfaces through the volume-preserving harmonic mean curvature flow (see Section 3). While \cite{Guilisun} treated the three-dimensional case, the higher-dimensional setting presents additional challenges since the curvature tensor is no longer completely determined by the Ricci tensor alone. To address this issue, we carefully decompose the curvature tensor into 3 components, including its scalar curvature, traceless Ricci curvature, and Weyl tensor, demonstrating that under our metric decay assumptions, the full curvature tensor can indeed be controlled by the Ricci curvature.

The uniqueness of the foliation constitutes a crucial aspect of our investigation. In contrast to the constant mean curvature case studied by Huisken-Yau and Qing-Tian, the CHMC case requires more delicate treatment. The principal difficulty stems from the fact that the linearized harmonic mean curvature operator $L$ (defined in \eqref{e-L}) fails to be self-adjoint, rendering standard elliptic theory inapplicable. To overcome this obstacle, we develop a two-fold strategy: first, we approximate $L$ by a scaled Laplacian operator; second, we establish sharp estimates for the symmetrization $S$ of $L$. The decay condition \eqref{metr} with $\delta>0$ proves essential here, as it ensures that the difference between $L$ and $S$ becomes negligible in our analysis. This explains why the local uniqueness result requires the stronger decay assumption $\delta>0$. Moreover, since the second variation of the harmonic mean curvature operator is of order $-3$, our uniqueness result is necessarily restricted to the three-dimensional case.

In a forthcoming paper, we will investigate the existence and uniqueness of constant harmonic mean curvature foliation using elliptic methods instead of parabolic approaches as in the CMC case \cite{Eich-Met}. This alternative framework is expected to yield stronger uniqueness results under more general conditions.

\vspace{.1in}

The paper is organized as follows. Section 2 presents preliminary materials and derives the evolution equations for geometric quantities under the volume-preserving harmonic mean curvature flow. In Section 3, we establish the long-time existence and exponential convergence of the flow. Section 4 addresses the construction of the foliation and the definition of the center of mass. Finally, Section 5 proves the uniqueness of the foliation for the three-dimensional case under the stronger decay assumption $\delta>0$.

\begin{acknowledge}
The first author is supported by the Initial Scientific Research Funding of No.X2450216. The second author is supported by Initial Scientific Research Fund of Young Teachers in Hefei university of Technology of No.JZ2024HGQA0122; The third author is supported by NSFC No. 12071352, 12271039.
\end{acknowledge}

\vspace{.2in}

\section{Preliminaries}

\vspace{.1in}

In this section, we establish fundamental estimates for curvature quantities and their derivatives under the metric decay conditions specified in \eqref{metr}.
Let us denote the Schwarzschild metric by
$ g_{\alpha\beta}^S = (1 + \frac{m}{2r^{n-2}})^{\f4{n-2}} \delta_{\alpha\beta}$. The following properties are well-known:

\begin{lem}\label{lem1.1}
     The Ricci curvature $R_{\alpha\beta}^S$ of $g_{\alpha\beta}^S$ is given by
$$R_{\alpha\beta}^S = \frac{(n-2)m}{r^n} \varphi^{4-2n} \left( \delta_{\alpha\beta} - n \frac{y^\alpha y^\beta}{r^2} \right),$$
where here and in the following $\varphi = \left(1 + \frac{m}{2r^{n-2}}\right)^{\frac{1}{n-2}}$.
\end{lem}
Thus one eigenvalue of the Ricci curvature (in radial direction) is given by $\f{(1-n)(n-2)m}{r^n\varphi^{2n-4}}$ and the other eigenvalues are $\f{m(n-2)}{r^n\varphi^{2n-4}}$. The scalar curvature equals zero.

\begin{lem}\label{lem1.2}
    Let $\bar{\nabla}$ and $\bar{R}_{\alpha\beta}$ denote the covariant differentiation and Ricci curvature respectively with respect to $\bar{g}$ given by \eqref{metr}. Then we have
$$\begin{aligned}
|\bar{R}_{\alpha\beta}-R_{\alpha\beta}^S|&\leq cc_0r^{-1-n-\delta},\\
|\bar{\nabla}_\gamma\bar{R}_{\alpha\beta}-\nabla_\gamma^SR_{\alpha\beta}^S|&\leq cc_0r^{-2-n-\delta},\\
|\bar{\nabla}_\rho\bar{\nabla}_\gamma\bar{R}_{\alpha\beta}-\nabla_\rho^S\nabla_\gamma^SR_{\alpha\beta}^S|&\leq cc_0r^{-3-n-\delta}.\end{aligned}$$
where c is an absolute constant.
\end{lem}

Consider an $(n-1)$-dimensional hypersurface $M$ embedded in $(N,\bar{g})$. Let $g=\{g_{ij}\}$ denote the induced metric on $M$, where Latin indices $i, j$ range from 1 to $(n-1)$. 
Let $\nu$ be the unit outward normal vector field along $M$ and denote the induced Levi-Civita connection on $M$ by $\nabla$. 
We often work in an orthonormal frame $e_1,\cdots,e_{n-1},\nu$ adapted to $M$, such that the second fundamental form $A=\{h_{ij}\}$ is given by 
$$
h_{ij}=-\langle\bar{\nabla}_{e_i}e_j,\nu\rangle=\langle\bar{\nabla}_{e_i}\nu,e_j\rangle.$$
Let $\lambda_1,\cdots,\lambda_{n-1}$ be the principal curvatures on $M$ and write 
$$H=g^{ij}h_{ij}=\sum_{i=1}^{n-1}\lambda_i,\quad|A|^2=h^{ij}h_{ij}=\sum_{i=1}^{n-1}\lambda_i^2$$
for the mean curvature and the square of the norm of the second fundamental form on $M$ respectively. 
Furthermore, let $\mathring{A}$ be the traceless second fundamental form such that
$$\mathring{h}_{ij}=h_{ij}-\frac1{n-1}Hg_{ij},$$
$$|\mathring{A}|^2=|A|^2-\frac1{n-1}H^2=\frac1{n-1}\sum_{i< j}(\lambda_i-\lambda_j)^2.$$
Moreover, the harmonic mean curvature is defined as
\[ F=\frac{1}{\sum_{i=1}^{n-1}\lambda_i^{-1}}. \]
The curvature of $N$ together with its derivatives and the second fundamental form of $M$ are related by
the equations of Gauss and Codazzi:
\[
\begin{aligned}
R_{ijkl}&=\bar{R}_{ijkl}+h_{il}h_{jk}-h_{ik}h_{jl},\\
-\bar{R}_{\nu ijk}&=\nabla_{k}h_{ij}-\nabla_{j}h_{ik}.
\end{aligned}
\]
 The derivatives of the curvature are related by 
\[
\bar{\n}_l\bar{R}_{\nu ijk}=\n_l\bar{R}_{\nu ijk}+h_{jl}\bar{R}_{\nu i\nu k}+h_{kl}\bar{R}_{\nu ij\nu}-h_{ml}\bar{R}_{mijk},
\]
and the commutation of derivatives are
\[
h_{ij,kl}=h_{ij,lk}+h_{im}R_{mjlk}+h_{jm}R_{milk}.
\]
We should state another consequence of the Codazzi equation, which will be utilized in dealing with the curvature estimates, especially to handle the extra terms coming from the derivative of $F$.
\begin{lem}\label{lem1.4}
If $w_{i}=-\bar{R}_{\nu lil}$ denotes the projection of $\bar{R}ic
(\nu,\cdot)$ onto $M$, then we have, for any $\eta>0$, the inequality
\[
|\nabla A|^{2}\geq\left(\frac{3}{n+1}-\eta\right)|\nabla H|^{2}-\frac{2}{n+1}\left(\frac{2}{(n+1)\eta}-\frac{n-1}{n-2}\right)|w|^{2}.
\]
Furthermore, we have the higher order version, for $k=1, 2, 3,\cdots$,
\[
|\n^kA|^2\geq\left(\frac{3}{n+1}-\eta\right)|\n^kH|^2-\frac{2}{n+1}\left(\frac{2}{(n+1)\eta}-\frac{n-1}{n-2}\right)|\n^{k-1}w|^2.
\]
\end{lem}
\begin{proof}
    The first result comes from \cite{Huisken1986ContractingCH}. To show the second statement, we first consider $k=2$ and then utilize induction argument. Similar to \cite{Huisken1986ContractingCH}, we define 
    \[
    \begin{aligned}
E_{ijkl}=&\frac{1}{n+1}(H_{il}g_{jk}+H_{jl}g_{ik}+H_{kl}g_{ij})-\frac{2}{(n+1)(n-2)}w_{kl}g_{ij}\\
&+\frac{n-1}{(n+1)(n-2)}(w_{il}g_{jk}+w_{jl}g_{ik}),
\end{aligned}
\]
\[F_{ijkl}=h_{ij,kl}-E_{ijkl}.
     \]
    We find that 
   \[ \begin{aligned}
    % \abs{E}^2=\f34H_{kl}^2-w_{kl}H_{kl}+w_{kl}^2,\\
  & \<E_{ijkl},h_{ij,kl}\>=\abs{E}^2\\
  =& \frac{3}{n+1}|\n^2H|^2+\frac{4}{n+1}H_{kl}w_{kl}+\frac{2(n-1)}{(n+1)(n-2)}|\n w|^2.
    \end{aligned}
      \]
      Thus, we conclude that $\<E_{ijkl},F_{ijkl}\>=0$. It follows that 
      \begin{align*}
      |\n^2A|^2 
      & = \abs{E}^2+\abs{F}^2\geq\abs{E}^2\\
      & \geq \left(\frac{3}{n+1}-\eta\right)|\n^2H|^2-\frac{4}{\eta(n+1)^2}|\n w|^2+\frac{2(n-1)}{(n+1)(n-2)}|\n w|^2\\
      & =\left(\frac{3}{n+1}-\eta\right)|\n^2H|^2-\frac{2}{n+1}\left(\frac{2}{\eta(n+1)}-\frac{n-1}{n-2}\right)|\n w|^2.
      \end{align*}
\end{proof}
The following Lemma is treated in \cite{Huisken1996DefinitionOC}, which will be used later.
\begin{lem}[Lemma 1.4 in \cite{Huisken1996DefinitionOC}]\label{phi}
Let M be a hypersurface in $(N,\bar{g})$ with principal curvatures $\lambda_i, i=1, 2, \cdots, n-1$. If $\hat{g}_{\alpha\beta}=\psi^2\bar{g}_{\alpha\beta}$ describes a conformal change of the ambient metric, then the eigendirections of the second fundamental form of M remain fixed and the new principal curvatures $\tilde{\lambda}_i$ are given by
$$\tilde{\lambda}_i=\psi^{-1}\lambda_i+\psi^{-2}\partial_\nu\psi.$$
In particular, the difference between the eigenvalues is conformally invariant:
$$(\tilde{\lambda}_i-\tilde{\lambda}_j)=\psi^{-1}(\lambda_i-\lambda_j).$$
\end{lem}
% \begin{proof}
%     The relation $\hat{g}_{\alpha\beta}=\psi^2\bar{g}_{\alpha\beta}$ implies that the Christoffel symbol $\hat{\Gamma}$ is given

% $$\hat{\Gamma}_{\beta\gamma}^\alpha=\Gamma_{\beta\gamma}^\alpha+\psi^{-1}\left(\delta_\beta^\alpha\partial_\gamma\psi+\delta_\gamma^\alpha\partial_\beta\psi-\bar{g}_{\beta\gamma}\bar{g}^{\alpha\rho}\partial_\rho\psi\right)$$
% and the result is then an immediate consequence of the definition of the second fundamental form above.
% \end{proof}
The proposition illustrate that near infinity the curvature of the more general metric $\bar{g}$ will be very close to the curvature of $g^S$.
\begin{prop}\label{prop2.1}
Let $M\hookrightarrow(N,\bar{g})$ be a hypersurface such that $\,r(y)\geq\f1{10}\max_{M}r=\f1{10}r_1, \forall y\in M.$ If for some constants $k_1, k_2$, there holds
\[|\mathring{A}|\leq k_1r_1^{-n-\delta},\quad\quad |\nabla\mathring{A}|\leq k_2r_1^{-n-\delta-1}.\]
Then there exists an absolute constant $ c>0$, such that the second fundamental form of $M$ with respect to the Euclidean metric statisfies
$$|\mathring{A}^e|\leq c(k_{1}+c_0)r_1^{-n-\delta},\quad\quad |\nabla^{e}A^e|\leq c(k_{2}+c_0)r_1^{-n-1-\delta},$$
provided that $r_1\geq c(c_0+k_{1})$. Moreover, $\exists\ r_{0}\in\R$ and a vector $\vec{a}\in\mathbb{R}^{n}$ such that
$$
|\lambda_i^e-r_{0}^{-1}|\leq c(k_{1}+k_{2}+c_{0})r_1^{-n-\delta},
$$
$$\left|{\nu}_{e}-\frac{y-\vec{a}}{r_0}\right|\leq c(k_{1}+k_{2}+c_{0})r_1^{-n-\delta+1},
$$
$$
|y-\vec{a}-r_{0}{\nu}_{e}|\leq c(k_{1}+k_{2}+c_{0})r_1^{-n-\delta+2},
$$
where $y$ and $\nu_e$ are the position vector and the unit normal of $M$ in $\mathbb{R}^n$ respectively.
\end{prop}

\begin{proof}
Since the Christoffel symbols of the metric \eqref{metr} and the Schwarzschild metric $g^S$ are close, the second fundamental form $A^S$ of $M$ with respect to the Schwarzschild metric $g^S$ satisfies
\[ |\mathring{A}^S-\mathring{A}|\leq cc_0r_1^{-n-\delta}, \quad |\n \mathring{A}^S-\n \mathring{A}|\leq cc_0r_1^{-n-1-\delta}, \]
which imply that
\[ |\mathring{A}^S|\leq c(c_0+k_1)r_1^{-n-\delta}, \quad |\n^S\mathring{A}^S|\leq c(c_0+k_2)r_1^{-n-\delta-1}, \]
provided that $r_1\geq c(c_0+k_1)$.
Let $\phi=\left(1+\frac{m}{2r^{n-2}}\right)^{\frac{1}{n-2}}$ and $\psi=\phi^2$. 
Then $g^S=\phi^4\delta=\psi^2\delta$.
From Lemma \ref{phi},  we obtain
\[ |\mathring{A}^{e}|\leq\psi|\mathring{A}^S|\leq c|\mathring{A}^S|\leq  c(c_0+k_1)r_1^{-n-\delta}, \]
\[ |\nabla^e\mathring{A}^{e}|\leq c(c_0+k_2)r_1^{-n-\delta-1}. \]
Using Lemma \ref{lem1.4} in Euclidean space, we get
\[ |\n^eA^e|\leq c(c_0+k_2)r_1^{-n-\delta-1}. \]

Note that at the point of $M$ where $\max_Mr=r_1$ is attained, we have 
$\lambda_i^e\geq r_1^{-1}$.
Since $r_1\leq 10r$, we see that there exists $r_0\in\mathbb{R}$ such that
\[ |\lambda_i^e-r_0^{-1}|\leq c(k_1+k_2+c_0)r_1^{-n-\delta},\quad i=1, 2, \cdots, n-1. \]
By the Gauss-Weingarten formula in the Euclidean space, we have
\[ \left|\frac{\partial}{\partial x_i}(\nu_e-r_0^{-1}y)\right|=\left|(h_{ij}^eg^{jk}-r_0^{-1}\delta_i^k)\frac{\partial y}{\partial x_k}\right|\leq c(k_1+k_2+c_0)r_1^{-n-\delta}. \]
Then there is a vector $\vec{a}\in \mathbb{R}^n$ such that the last two estimates hold.
\end{proof}

Next we define a set of round hypersurfaces in $(N, \bar{g})$. 
For $\sigma\geq1$, $B_1, B_2, B_3>0$, we set  
{\small $$\mathcal{B}_{\sigma} (B_1, B_2, B_3) = \left\{ M \subset N \mid |r-\sigma| \leq  B_1, \; |\mathring{A}| \leq B_2 \sigma^{-n-\delta}, \; |\nabla \mathring{A}| \leq B_3 \sigma^{-n-1-\delta} \right\}.$$}

For $\sigma$ sufficiently large, we can use Proposition \ref{prop2.1} to estimate the principal curvatures of hypersurfaces in $\mathcal{B}_{\sigma} (B_1, B_2, B_3)$.

\begin{prop}\label{prop2.2}
    Let $M$ be a hypersurface in $\mathcal{B}_{\sigma} (B_1, B_2, B_3)$.
    Suppose that $\sigma\geq c(c_0+B_1+B_2)$ is such that Proposition \ref{prop2.1} holds. 
    Then the principal curvature $\lambda_i$, for $i=1, 2, \cdots, n-1$, of $M$ satisfies
    \[ \left|\lambda_i-\frac{1}{r_{0}}+\frac{(n-1)m}{n-2}r_{0}^{1-n}-mn\<\vec{a},\nu_e\>_er_0^{-n}-\frac{nm^2}{4(n-2)^{2}}r_0^{3-2n}\right|\leq C{\sigma}^{-n-\delta} \]
    and the mean curvature satisfies
    \[ \left|H-\frac{n-1}{r_{0}}+\frac{(n-1)^2m}{n-2}r_{0}^{1-n}-\f{mn(n-1)\<\vec{a},\nu_e\>_e}{r_0^n}-\frac{n(n-1)m^2}{4(n-2)^{2}}r_{0}^{3-2n}\right|\leq C\sigma^{-n-\delta}, \]
    provided that $\sigma\geq c(c_0+B_1^2+B_2)$, where $C=c(c_0^2+B_2+B_3)$.
\end{prop}

\begin{proof}
Let $$\varphi=\left(1+\frac{m}{2r^{n-2}}\right)^{\frac{1}{n-2}},\, \psi=\varphi^2,\, g^S=\varphi^{4}\delta=\psi^2\delta.
$$

From Lemma \ref{prop2.1}, we have
\begin{equation}\label{lam}
 \begin{gathered}
\lambda_i^e=\frac{1}{r_{0}}+\delta_{1},\quad |\delta_{1}|\leqslant c(k_{1}+k_2+c_{0})r_1^{-n-\delta},\\
y=\vec{a}+r_{0}\nu_{e}+\delta_{2}\nu_e,\quad
|\delta_{2}|\leq c(k_{1}+k_2+c_{0})r_1^{-n-\delta+2}.
\end{gathered} 
\end{equation}
By Lemma \ref{phi} and \eqref{lam}, we calculate directly that
\begin{align*}
\lambda_{i}^S&=\lambda_{i}^{e}\varphi^{-2}+\varphi^{-4}\partial_{\nu_e}\left(1+\frac{m}{2r^{n-2}}\right)^{\frac{2}{n-2}}\\
&=\varphi^{-2}\lambda_{i}^{e}+\varphi^{-4}\left(\frac{2}{n-2}\right)\left(1+\frac{m}{2r^{n-2}}\right)^{\frac{2}{n-2}-1}\cdot\frac{m}{2}(-n+2)r^{1-n}\<\f yr,{\nu}_{e}\>_{e}\\
&=\varphi^{-2}\lambda_i^e-m\varphi^{-n}r^{-n}\<y,\nu_e\>_e\\
&=\left(\frac{1}{r_{0}}+\delta_{1}\right)\varphi^{-2}-m\varphi^{-n}r^{-n}\<\vec{a}+r_{0}{\nu}_{e}+\delta_{2}{\nu}_{e}, {\nu}_{e}\>_e\\
&=\frac{1}{r_{0}}\varphi^{-2}-m\varphi^{-n}r^{-n}\<\vec{a}, {\nu}_{e}\>_e-m\varphi^{-n}r^{-n}r_{0}+\delta_{1}\varphi^{-2}-m\varphi^{-n}r^{-n}\delta_{2}.
\end{align*}
This implies
\[
\left|\lambda_i^S-\frac{1}{r_{0}}\varphi^{-2}+m\varphi^{-n}r^{-n}\<\vec{a}, {\nu}_{e}\>_e+m\varphi^{-n}r^{-n}r_{0}\right|\leq C\sigma^{-n-\delta}.
\]
Note that
\[
|\lambda_i^S-\lambda_i|\leq cc_0\sigma^{-n-\delta},\quad (1+x)^{\alpha}=1+\alpha x+\frac{\alpha(\alpha-1)}{2}x^{2}+O(x^{3}).
\]
It follows that 
$$\varphi^{-2}=\left(1+\frac{m}{2r^{n-2}}\right)^{\frac{-2}{n-2}}=1-\frac{m}{(n-2)r^{n-2}}+\frac{nm^{2}}{4(n-2)^{2}r^{2n-4}}+O(r^{-3n+6}),$$
and
$$
\varphi^{-n}=\left(1+\frac{m}{2r^{n-2}}\right)^{\frac{-n}{n-2}}=1-\frac{nm}{2(n-2)r^{n-2}}+O(r^{-2n+4}),
$$
which implies 
$$\left|\lambda_i-\frac{1}{r_{0}}+\frac{m}{(n-2)r_{0}r^{n-2}}-\frac{nm^{2}}{4(n-2)^{2}r_{0}r^{2n-4}}+\f m{r^n}\<\vec{a},\nu_e\>_e+\f {mr_0}{r^n}\right|\leq c\sigma^{-n-\delta}.
$$

By the Cosine Rule,  we see that $r=r_0+\<\vec{a},\nu_e\>_e+\delta_3$ for $\abs{\delta_3}\leq\f{\abs{\vec{a}}}{2r_0}$. Then we have
\begin{align*}
&r^{2-n}=r_{0}^{2-n}\left(1+\frac{\<\vec{a},\nu_e\>_e+\delta_3}{r_0}\right)^{2-n}\\
=&r_{0}^{2-n}\left(1+(2-n)\frac{\<\vec{a},\nu_e\>_e}{r_0}+O(r_0^{-2})\right)\\
=&r_{0}^{2-n}+(2-n)r_0^{1-n}\<\vec{a},\nu_e\>_e+O(r_0^{-n}),
\end{align*}
and 
\[ r^{-n}=r_{0}^{-n}-n\<\vec{a},\nu_e\>_er_0^{-n-1}+O(r_0^{-n-2}). \]
Therefore, we obtain that
\begin{align*}
\left|\lambda_i-\frac{1}{r_{0}}+\frac{(n-1)m}{n-2}r_{0}^{1-n}-mn\<\vec{a},\nu_e\>_er_0^{-n}-\frac{nm^2}{4(n-2)^{2}}r_0^{3-2n}\right|\leq C\sigma^{-n-\delta}.
\end{align*}
and consequently,
{\small \[ \left|H-\frac{n-1}{r_{0}}+\frac{(n-1)^2m}{n-2}r_{0}^{1-n}-\f{mn(n-1)\<\vec{a},\nu_e\>_e}{r_0^n}-\frac{n(n-1)m^2}{4(n-2)^{2}}r_{0}^{3-2n}\right|\leq C\sigma^{-n-\delta}. \]}
\end{proof}

\vspace{.1in}

\section{Existence of constant harmonic mean curvature hypersurfaces}

\vspace{.1in}

In this section, we establish the existence of constant harmonic mean curvature hypersurfaces in asymptotically Schwarzschild spaces via the volume-preserving harmonic mean curvature flow, with initial data given by coordinate spheres.
Precisely, for each $\sigma>0$, let $S_\sigma(0)$ denote the coordinate sphere of radius $\sigma$ centered at the origin in $(N,\bar{g})$ given by a map
$$
\phi_0^\sigma:S^{n-1}\to N,\quad\phi_0^\sigma(S^{n-1})=S_\sigma(0).
$$
We aim to find a one-parameter family of maps $\phi_t^\sigma=\phi^\sigma(\cdot,t)$ solving the initial value problem
\bel{flow}
\begin{cases}
\frac{d}{dt}\phi^{\sigma}(p,t)&=(f-F)\nu(p,t),\quad t\geq0,\:p\in S^{n-1},\\
\phi^{\sigma}(0)&=\phi_{0}^{\sigma},
\end{cases}
\qe
where $F(p, t)=\frac{1}{\sum_{i=1}^{n-1}\lambda_i^{-1}}$ and $\nu(p, t)$ represent the harmonic mean curvature and the outward unit normal vector of the surface $\Sigma_t=\phi_t^\sigma(S^{n-1})$ respectively. Here we denote $|\Sigma_t|$ and $f(t)=\frac{\int_{\Sigma_t}Fd\mu_t}{|\Sigma_t|}$ as the area and the average of the harmonic mean curvature. 

\begin{thm}\label{main}
     If $(N, \bar{g})$ is the asymptotically Schwarzschild manifold given by \eqref{metr} with $\delta\geq 0$, then there is $\sigma_0$ depending only on $c_0$ such that for all $\sigma\geq\sigma_0$, the volume preserving harmonic mean curvature flow \eqref{flow} has a smooth solution for all times $t\geq0$.
    As $t\rightarrow\infty$, the hypersurfaces $\Sigma_t$ converge exponentially fast to a hypersurface $\Sigma_{\sigma}$ of constant harmonic mean curvature $f_{\sigma}$.
\end{thm}

To show the long time existence and exponential convergence of the flow \eqref{flow}, we will first prove that $\Sigma_t$ remains in $\mathcal{B}_{\sigma} (B_1, B_2, B_3)$ for $\sigma\geq \sigma_0$.
Similar to the proof of Huisken-Yau in \cite{Huisken1996DefinitionOC}, there hold the following propositions, where we control the position of $\Sigma_t$ during the evolution by estimates which are independent of $B_1$.

\begin{prop}\label{3.4}(Proposition 3.4 in \cite{Huisken1996DefinitionOC})
    Suppose that $\Sigma_t$ is a smooth solution of the flow \eqref{flow} contained in $\mathcal{B}_{\sigma}(B_1, B_2, B_3)$ for all $t\in [0, T]$. Assume that $\sigma\geq c(C_0+B_1+B_2)$ is such that Proposition \ref{prop2.1} applies and $r_0(t)$ is defined as in that result. Then there is an absolute constant $c$ such that
    \[|r_0(t)-\sigma|\leq c(c_0+B_2+B_3) \]
    holds uniformly in $[0, T]$, provided that $\sigma\geq c(c_0+B_1+B_2)$.
\end{prop}
The proof is similar as \cite{Huisken1996DefinitionOC}, since the flow we are investigating is volume preserving, and since the metric is asymptotically flat, thus the volume is almost Euclidean, we then conclude the radius remains invariant up to constants.
The following Proposition shows that the evolution hypersurfaces will remain in $\mathcal{B}_\sigma(B_1,B_2,B_3)$, and thus will not drift away.
\begin{prop}\label{3.5}
    Suppose that the solution $\Sigma_t$ of \eqref{flow} is contained in $\mathcal{B}_{\sigma}(B_1, B_2, B_3)$ for all $t\in [0, T]$. Then there is an absolute constant $c$ such that
    \[ \max_{\Sigma_t}r\leq \sigma+c(m^{-1}+1)(c_0^2+B_2+B_3) \]
     holds uniformly in $[0, T]$, provided that $\sigma\geq c(c_0^2+B_1^2+B_2+B_3)$.
\end{prop}

\begin{proof}
    Let $D>0$ and assume $\max_{\Sigma_t}r<\sigma+D$ is violated for the first time at $(y_0, t_0),\  t_0>0$.
    At that point $r(y_0, t_0)=\max_{\Sigma_t}r=\sigma+D$, $\langle y_0, \nu_e\rangle_e=r$ and 
   $(f-F)\langle \nu, \nu_e\rangle_e\geq0$.
    From Proposition \ref{prop2.1}, it follows that $\langle \nu, \nu_e\rangle_e\geq\frac{1}{2}$ and $\langle \Vec{a}, \nu_e\rangle_e\geq\frac{1}{2}|\vec{a}|$ at $(y_0, t_0)$ provided $\sigma\geq c(C_0+B_2+B_3)$. 
    So we have $F\leq f$ at $(y_0, t_0)$, but we get from Proposition \ref{prop2.2} that
    \begin{equation*}
    \begin{split}
       & F = \frac{\lambda_1 \cdots \lambda_{n-1}}{\sum_{1\leqslant i_1<\cdots<i_{n-2}\leqslant {n-1}} \lambda_{i_1} \cdots \lambda_{i_{n-2}}} \\
       =& \frac{1}{(n-1) r_0} - \frac{m}{(n-2)} r_0^{1-n} + \frac{mn \langle \vec{a}, \nu_e \rangle_e}{n-1} r_0^{-n} + \frac{nm^2 r_0^{3-2n}}{4(n-1)(n-2)^2} + O(\sigma^{-n-\delta}).
    \end{split}
\end{equation*}
    Thus, we get at $(y_0,t_0)$,
    \begin{equation}\label{3.2}
    \begin{split}
        F-f =& \frac{nm}{(n-1)r_0^n}\left(\langle \vec{a}, \nu_e\rangle_e-\oint_{\Sigma_t}\langle \vec{a}, \nu_e\rangle_ed\mu_t\right)+O(\sigma^{-n-\delta})\\
        \geq& \frac{nm|\vec{a}|}{2(n-1)}\sigma^{-n}-c(c_0^2+B_2+B_3)\sigma^{-n-\delta},
    \end{split}
    \end{equation} 
    provided that $\sigma\geq c(c_0^2+B_1^2+B_2+B_3)$. At $(y_0, t_0)$, in view of Proposition \ref{3.4},
     we have
     \[ |\vec{a}|\geq D-c(c_0+B_2+B_3). \]
     Combining with \eqref{3.2}, we finally get
     \[ 0\geq F-f\geq \frac{nmD}{2(n-1)}\sigma^{-n}-c(c_0^2+B_2+B_3)\sigma^{-n}, \]
     which yields a contradiction if $D>c(m^{-1}+1)(c_0^2+B_2+B_3)$.
\end{proof}

We are now in a position to give the curvature estimates in order to derive the long time existence. We denote the linearized operator by $\mathcal{L}f=F^{kl}f_{kl}$ as in \cite{Guilisun}. 
We can calculate the evolution equations directly.
\begin{align*}
(\frac{\pa}{\pa t}-\mathcal{L})|A|^2=&-2F^{kl}h_{ij,k}h_{ij,l}-2ftrA^3+2F^{kl}h_{mk}h_{ml}|A|^{2}+2fh_{ij}\bar{R}_{\nu i\nu j}\\
&-2|A|^{2}F^{kl}\bar{R}_{\nu k\nu l}-4F^{kl}h_{ij}\big(h_{jm}\bar{R}_{mkli}-h_{km}\bar{R}_{mjil}\big)\\
&-2F^{kl}h_{i j}(\bar{\n}_l\bar{R}_{\nu jk i}+\bar{\n}_i\bar{R}_{\nu kl j})+2F^{kl,pq}h_{ij}\nabla_ih_{kl}\nabla_{j}h_{pq},
\end{align*}
\begin{align*}
(\frac{\pa}{\pa t}-\mathcal{L})H^{2}=&2H(\frac{\pa}{\pa t}-\mathcal{L})H-2F^{kl}\nabla_{k}H\nabla_lH\\
=&2HF^{kl,pq}\nabla_ih_{pq}\nabla^{i}h_{kl}-2F^{pq}H_pH_q-2Hf(|A|^{2}+\bar{R}ic(\nu,\nu))\\
&-4H F^{kl}(h_{m}^{i}\bar{R}_{mkl i}-h_{lm}\bar{R}_{m iik})+2H^{2}F^{kl}(h_{mk}h_{ml}-\bar{R}_{\nu k\nu l})\\
&-2HF^{kl}(\bar{\nabla}_l\bar{R}_{\nu ik i}+\bar{\nabla}^j\bar{R}_{\nu klj}),
\end{align*}
which implies 
\begin{equation}\label{a0}
\begin{aligned}
    &(\frac{\partial}{\partial t}-\mathcal{L})|\mathring{A}|^2\\=& -2F^{kl}\mathring{h}_{ij,k}\mathring{h}_{ij,l}-\frac{2f}{n-1}\sum_{i<j}(\lambda_i+\lambda_j)(\lambda_i-\lambda_j)^2+2F^{kl}h_{mk}h_{ml}|\mathring{A}|^2\\
    &+2F^{kl,pq}\mathring{h}_{ij}\nabla_ih_{kl}\nabla_jh_{pq}
    +2f\mathring{h}_{ij}\bar{R}_{\nu i\nu j}-2F^{kl}|\mathring{A}|^2\bar{R}_{\nu k\nu l}\\
    &-2\mathring{h}_{ij}F^{kl}\left(\bar{\n}_l\bar{R}_{\nu jki}+\bar{\n}_i\bar{R}_{\nu kl j}\right)
    -4F^{kl}\mathring{h}_{ij}\left(h_{jm}\bar{R}_{mkli}-h_{km}\tilde{R}_{mjil}\right).
\end{aligned}
\end{equation}
We define $e=\f{|\mathring{A}|^2}{H^2}=\f{A^2}{H^2}-\f1{n-1}$ and calculate as
\begin{align*}
(\frac{\partial}{\partial t}-\mathcal{L})e=&\frac{2}{H}F^{lk}\nabla_{l}H\n_ke-\frac{2}{H^4}|\nabla_{i}Hh_{kl}-H\nabla_{i}h_{kl}|_{F}^{2}+\frac{2f}{H}\bar{R}ic(\nu,\nu)e\\
&+\underbrace{\frac{2f}{H^{3}}(\abs{A}^{4}-Htr(A^3))}_{(2)}
-\underbrace{\frac{4F^{lk}}{H^{2}}(h_{ij}-\frac{\abs{A}^2}{H}g_{ij})(h_{jm}\bar{R}_{mkli}-h_{km}\bar{R}_{mjil})}_{(4)}\\
&+\underbrace{\frac{2f}{H^{2}}\mathring{h}_{ij}\bar{R}_{\nu i\nu j}}_{(3)}-\underbrace{\frac{2F^{kl}}{H^{2}}(h_{ij}-\frac{\abs{A}^2}{H}g_{ij})(\bar{\n}_l\bar{R}_{\nu jki}+\bar{\n}_{i}\bar{R}_{\nu klj})}_{(5)}\\
&+\underbrace{\frac{2}{H^{2}}(h_{ij}-\frac{\abs{A}^2}{H}g_{ij})(F^{kl,pq}\n_ih_{kl}\n_jh_{pq})}_{(6)},
\end{align*}
where 
\[
|\nabla_{i}Hh_{kl}-\nabla_{i}h_{kl}H|_F^2=F^{kl}(h_{ij,k}h_{ij,l}H^2+H_kH_l\abs{A}^2-2h_{ij,k}h_{ij}H_lH).
\]
To derive the desired curvature estimates, we need the following lemma.

\begin{lem}\label{curv}
 Let $\Sigma_t$ be a smooth solution of \eqref{flow} contained in $\mathcal{B}_\sigma(B_1,B_2,B_3)$,
for $t\in[0,T]$. Then there is an absolute constant c such that for $\sigma\geq c (c_{0}^{2}+B_1^2+B_2)$
{\renewcommand{\theequation}{\roman{equation}}\setcounter{equation}{0}
\begin{gather}
 |\lambda_i - \frac{1}{\sigma}| \leq C\sigma^{1-n}, \quad \quad |F - f| \leq C\sigma^{1-n},\\
  \frac{2f}{H^3} (|{A}|^4 - HtrA^3 ) \leq -\frac{1}{(n-1)^3} |\mathring{A}|^2,\\
 \left| \frac{f}{H^2} \mathring{h}_{ij} \bar{R}_{\nu i\nu j} \right| \leq C \sigma^{-n-\delta} |\mathring{A}|,\\
\left|\frac{1}{H^2}F^{kl}(h_{ij}-\frac{|A|^2}{H}g_{ij})(h_{jm}\bar{R}_{mkli}-h_{lm}\bar{R}_{mijk})\right|\leq c|\mathring{A}|\sigma^{2-2n-\delta},\\
\left|\frac{F^{kl}}{H^{2}}\left(h_{ij}-\frac{|A|^{2}}{H}g_{ij}\right)\left(\bar{\nabla}_{l}\bar{R}_{\nu jki}+\bar{\nabla}_{i}\bar{R}_{\nu klj}\right)\right|\leq C\sigma^{-n-\delta}|\mathring{A}|,\\
\left|\frac{2}{H^2}\left(h_{ij}-\frac{|A|^2}{H}g_{ij}\right)F^{kl,pq}\nabla_ih_{kl}\nabla_{j}h_{pq}\right|\leq C|\mathring{A}|\sigma^{-2n+1-2\delta}.
\end{gather}\setcounter{equation}{5}}
\end{lem}

\begin{proof}
The first inequality follows from Proposition \ref{prop2.2}.
The second inequality is obtained by 
  $$\frac{2f}{H^3} (|{A}|^4 - HtrA^3 )= -\frac{2f}{H^3} \sum_{i<j} \lambda_i \lambda_j (\lambda_i - \lambda_j)^2 \leq -\frac{1}{(n-1)^3} |\mathring{A}|^2.$$
Since the Schwarzschild metric $g^S$ has zero scalar curvature and is conformally flat, the scalar and Weyl curvature of $\bar{g}$ are of order $O(\sigma^{-n-1-\delta})$. By the curvature decomposition, we obtain that
$$  \bar{R}_{\nu i\nu j} = \frac{(\bar{R}ic \odot \bar{g})_{\nu i\nu j}}{n-2} + O(\sigma^{-1-n-\delta})=\frac{-\bar{R}_{ij} - \bar{R}_{\nu\nu} g_{ij}}{n-2} + O(\sigma^{-1-n-\delta}),$$
and consequently, 
$$\bar{R}_{\nu i\nu j} \mathring{h}_{ij} = -\frac{\bar{R}_{ij} \mathring{h}_{ij}}{n-2} + O(\sigma^{-1-n-\delta}) |\mathring{A}|.$$
Note that
\[
\begin{split}
   &\bar{R}_{ij} \mathring{h}_{ij}=\sum_{i=1}^{n-1}\bar{R}_{ii}(\lambda_i-\frac{H}{n-1})= \frac{1}{n-1}\sum_{i=1}^{n-1}\bar{R}_{ii}\sum_{j=1}^{n-1}(\lambda_i-\lambda_j)\\
   =& \frac{1}{n-1} \left[ \sum_{i<j} \bar{R}_{ii} (\lambda_i - \lambda_j) + \sum_{i>j} \bar{R}_{ii} (\lambda_i - \lambda_j) \right] \\
   =& \f1{n-1}\sum_{i<j}\left( \bar{R}_{ii} - \bar{R}_{jj} \right) (\lambda_i - \lambda_j).
\end{split}
\]
Since
$$
|\bar{R}ic - \bar{R}ic^S| \leq C \sigma^{-1-n-\delta}, $$ we get
$$ |\bar{R}_{ij} \mathring{h}_{ij}| \leq C \sigma^{-1-n-\delta} |\mathring{A}|,$$
which implies the third inequality.
To show the fourth inequality, we will also decompose the curvatures as following
$$\bar{R}_{mki l}=\frac{1}{n-2}[\bar{R}_{ml}\bar{g}_{ki}+\bar{R}_{ki}\bar{g}_{ml}-\bar{R}_{mi}g_{kl}-\bar{R}_{kl}g_{mi}]+O(\sigma^{-1-n-\delta}).$$
Then
$$
h_{jm}\bar{R}_{mki l}=\frac{1}{n-2}[h_{jm}\bar{R}_{ml}\bar{g}_{ki}+h_{jl}\bar{R}_{ki}-h_{jm}\bar{R}_{mi}\bar{g}_{kl}-h_{ij}\bar{R}_{kl}]+O(\sigma^{-1-n-\delta}),$$
$$h_{lm}\bar{R}_{mijk}=\frac{1}{n-2}[h_{lm}\bar{R}_{mk}\bar{g}_{ij}+h_{lk}\bar{R}_{ij}-h_{lm}\bar{R}_{mj}\bar{g}_{ik}-h_{lj}\bar{R}_{ik}]+O(\sigma^{-1-n-\delta}).$$
Note that 
$$
h_{ij}-\frac{|A|^2}{H}g_{ij}=\mathring{h}_{ij}+\left(\frac{H}{n-1}-\frac{|\mathring{A}|^2+\f{H^2}{n-1}}{H}\right)g_{ij}=\mathring{h}_{ij}-\frac{|\mathring{A}|^2}{H}g_{ij},$$
thus, we have
\begin{equation*}
    \begin{split}
       & F^{kl}(h_{ij}-\frac{|A|^2}{H}g_{ij})(h_{jm}\bar{R}_{mkil}+h_{lm}\bar{R}_{mijk})\\
       ={}&F^{kl}(\mathring{h}_{ij}-\frac{|\mathring{A}|^2}{H}g_{ij})\frac{1}{n-2}[h_{jm}\bar{R}_{ml}\bar{g}_{ki}+h_{lm}\bar{R}_{mk}\bar{g}_{ij}+h_{lk}\bar{R}_{ij}-h_{ij}\bar{R}_{kl}\\
       &-h_{jm}\bar{R}_{mi}\bar{g}_{kl}-h_{lm}\bar{R}_{mj}\bar{g}_{ik}+O(\sigma^{-1-n-\delta})]\\
={}&\f{F^{kl}}{n-2}(\mathring{h}_{ij}-\frac{|\mathring{A}|^2}{H}g_{ij})\big[\mathring{h}_{jm}\bar{R}_{ml}\bar{g}_{ki}+\mathring{h}_{lm}\bar{R}_{mk}\bar{g}_{ij}+\mathring{h}_{kl}\bar{R}_{ij}-\mathring{h}_{ij}\bar{R}_{kl}\\
&-\mathring{h}_{jm}\bar{R}_{mi}g_{kl}-\mathring{h}_{lm}\bar{R}_{mj}g_{ik}+O(\sigma^{-1-n-\delta})\big].
\end{split}
\end{equation*}
Therefore, there holds
\[
\left|\frac{1}{H^2}{F}^{kl}(h_{ij}-\frac{|A|^2}{H}g_{ij})(h_{jm}\bar{R}_{mkli}-h_{lm}\bar{R}_{mijk})\right|\leq c|\mathring{A}|^2\sigma^{2-n}\leq c|\mathring{A}|\sigma^{2-2n-\delta}.
   \]

Let $ w_{i}=\bar{R}_{\nu i}$, then $\abs{w}\leq C\sigma^{-n-1-\delta}$.
And
\begin{align*}
\bar{R}_{\nu jki}&=\frac{(\bar{R}ic\odot g)_{\nu jki}}{n-2}+O\left(\sigma^{-n-1-\delta}\right)\\
&=\frac{1}{n-2}\left[\bar{R}_{\nu i}g_{jk}-\bar{R}_{\nu k}g_{ij}\right]+O\left(\sigma^{-n-1-\delta}\right)=O(\sigma^{-n-1-\delta}),
\end{align*}
gives the fifth inequality
\[\left|\frac{F^{kl}}{H^{2}}\left(h_{ij}-\frac{|A|^{2}}{H}g_{ij}\right)\left(\bar{\nabla}_{l}\bar{R}_{\nu jki}+\bar{\nabla}_{i}\bar{R}_{\nu klj}\right)\right|\leq C\sigma^{-n-\delta}|\mathring{A}|.
\]
By Lemma \ref{lem1.4}, we have
\begin{gather*}
|\nabla A|^2\leq C_1|\n\mathring{A}|^2+C_2|w|^{2}\leq C\sigma^{-2n-2-2\delta}.
\end{gather*}
Then the final inequality follows from Proposition \ref{prop2.2}.
\end{proof}

 \begin{prop}\label{prop1.8}
Suppose that the solution $\Sigma_t$ is contained in $\mathcal{B}_\sigma(B_1,B_2,B_3)$ for $t\in[0,T]$. Then there is an absolute constant c such that for $\sigma\geq c(c_0^2+B_1^2+B_2+B_3)$ we have the estimate
$$
|\mathring{A}|^2\leq cc_0^2\sigma^{-2n-2\delta}
$$
everywhere in $[0,T]$.
\end{prop}
\begin{proof}
From Lemma \ref{curv}, we have
    \[
     \left(\frac\pa{\pa t}-\mathcal{L}\right)e\leq\frac{2}{H}F^{kl}\nabla_{l}H\nabla_{k}e-\frac{1}{(n-1)^3}|\mathring{A}|^{2}+c|\mathring{A}|\sigma^{-n-\delta},
     \]
     where we use the facts that $F^{ij}$ is uniformly elliptic and $\bar{R}ic(\nu, \nu)$ is negative.
     Then the maximum principle implies $|\mathring{A}|\leq C\sigma^{-n-\delta}$.
\end{proof}
   In order to get higher order control of derivatives of $\mathring{A}$, we need the following lemma, which gives the control of usual derivatives of harmonic mean curvature. The proof uses the facts that the hypersurfaces remains in $\mathcal{B}_\sigma(B_1,B_2,B_3)$, and hence the principal curvatures together with their derivatives can be controlled. 
\begin{lem}[Lemma 3.9 in \cite{Guilisun}]\label{lem1.9}
    Let $\Sigma_t$ be a smooth solution of \eqref{flow} contained in $\mathcal{B}_\sigma(B_1, B_2, B_3)$ for $t\in[0, T]$. Then there is an absolute constant $c$ such that  
    \[
    |d^kF|\leq c c_0^2\sigma^{k-1},
    \]
    for $\sigma\geq c(c_0^2+B^2_1+B_2)$. Here $d^kF$ is the k-th derivative with respect to the second fundamental form. 
\end{lem}
We can now derive the evolution equation for the first derivative of the second fundamental form by direct calculations as follows

\begin{lem}\label{na}
    If $\Sigma_t$ is a solution of \eqref{flow} contained in $\mathcal{B}_\sigma(B_1, B_2,B_3)$, then
there is an absolute constant $c_1$ such that
\[
\pd t|\n\mathring{A}|^2 \leq\mathcal{L}|\n\mathring{A}|^2+c_1|\nabla\mathring{A}\mid^2\sigma^{-2}+c_1c_0|\nabla\mathring{A}\mid\sigma^{-n-3-\delta},
\]
provided $\sigma\geq c_1 (C_0^2+B_1^2+B_2)$.
\end{lem}
We are now in a position to give the derivative estimate of $\mathring{A}$.
\begin{prop}\label{3.10}
    Suppose that the solution $\Sigma_t$ of \eqref{flow} is contained in $\mathcal{B}_\sigma(B_1,B_2, B_3)$ for $t\in[0,T]$. Then there is an absolute constant c such that for $\sigma \geq c( c_0^2+B_1^2+B_2+ B_3)$ the estimate
$$
|\nabla\mathring{A}|\leq cc_0\sigma^{-n-1-\delta}
$$
holds everywhere in $[0,T]$.
\end{prop}

\begin{proof}
    By \eqref{a0}, we see that
    \[ (\frac{\partial}{\partial t}-\mathcal{L})|\mathring{A}|^2\leq c\sigma^{-2n-2-2\delta}-\frac{1}{(n-1)^2}|\n \mathring{A}|^2. \]
    From Lemma \ref{na}, we have
    \[ (\frac{\partial}{\partial t}-\mathcal{L})(|\n \mathring{A}|^2+c_1|\mathring{A}|^2\sigma^{-2})\leq -c_3|\n \mathring{A}|^2\sigma^{-2}+c_4\sigma^{-2n-4-2\delta}. \]
    Using the maximum principle and Proposition \ref{prop1.8} gives the desired estimate. 
\end{proof}

Once we have obtained that the hypersurface $\Sigma_t$ is contained in $\mathcal{B}_\sigma(B_1,B_2, B_3)$ for all  $t$, higher derivative estimates and long-time existence are obtained. And hence the long time existence follows combining with  usual compactness argument, we have only to show that the solution hypersurfaces $\Sigma_t$ converge exponentially fast to a constant harmonic mean curvature limiting hypersurface.

\begin{lem}\label{mu1}
    Suppose $M$ is a hypersurface of $N$, then for every u with $\int_Mu=0$, we have
    \[
    \int_Mu\mathcal{L}u\leq-\left(\f1{(n-1)\sigma^2}-\f{2m}{(n-1)(n-2)\sigma^n}-c\sigma^{-n-1}\right)\int_Mu^2,
    \]
    and
    \[
    \int_MF^{ij}u_{i}u_j\geq\left(\f1{(n-1)\sigma^2}-\f{2m}{(n-1)(n-2)\sigma^n}-c\sigma^{-n-1}\right)\int_Mu^2.
    \]
\end{lem}
\begin{proof}
    By Proposition \ref{prop2.2}, $H=\f{n-1}\sigma-\f{(n-1)^2m}{n-2}\sigma^{1-n}+O(\sigma^{-n})$, which implies 
    \begin{align*}
        H^2=&\f{(n-1)^2}{\sigma^2}-\f{2(n-1)^3m}{n-2}\sigma^{-n}+O(\sigma^{-n-1})\\
        \geq&\f{(n-1)^2}{\sigma^2}-\f{2(n-1)^3m}{n-2}\sigma^{-n}-C\sigma^{-n-1}.
    \end{align*}
By Lemma 5.6 and Proposition 5.7 of \cite{Eich-Met}, we have that 
\begin{align*}
    \mu_{Lap}&\geq\f1{n-1}\left[\f{(n-1)^2}{\sigma^2}-\f{2(n-1)^3m}{n-2}\sigma^{-n}-C\sigma^{-n-1}\right]+\f{2(n-1)m}{\sigma^n}+O(\sigma^{-n-1})\\
    &\geq\f{n-1}{\sigma^2}-2\f{(n-1)m}{n-2}\sigma^{-n}-c\sigma^{-n-1}.
\end{align*}
Now we write that $\mathcal{L}u=F^{kl}u_{kl}=\f1{(n-1)^2}\Delta u+(F^{kl}-\f{\delta^{kl}}{(n-1)^2})u_{kl}$, which gives that
\begin{align*}
    \int_Mu\mathcal{L}u=&\f1{(n-1)^2}\int_Mu\Delta u+\int_Mu\left(F^{kl}-\f{\delta^{kl}}{(n-1)^2}\right)u_{kl}\\
    =&-\f1{(n-1)^2}\int_M\abs{\n u}^2-\int_Muu_kF^{kl}_{,l}-\int_M\left(F^{kl}-\f{\delta^{kl}}{(n-1)^2}\right)u_ku_l.
\end{align*}
Since $F=\f{\sigma_{n-1}}{\sigma_{n-2}}=\frac{1}{\sum_{i=1}^{n-1}\lambda_i^{-1}}$, thus 
\[
F^{ii}=\f{\pa F}{\pa\lambda_i}=\f{\lambda_i^{-2}}{(\sum_j\lambda_j^{-1})^2}.
\]
It then implies 
\begin{align*}
    F^{ii}-\f1{(n-1)^2}=&\f{(n-1)^2\lambda_i^{-2}-(\sum_j\lambda_j^{-1})^2}{(n-1)^2(\sum_j\lambda_j^{-1})^2}\\
    =&F^2\f{(n-1)^2\lambda_i^{-2}-\sum_{j,k}\lambda_j^{-1}\lambda_k^{-1}}{(n-1)^2}\\
    =&F^2\f{\sum_{j,k}(\lambda_i^{-2}-\lambda_j^{-1}\lambda_k^{-1})}{(n-1)^2}\\
    =&\f{F^2}{(n-1)^2}\sum_{j,k}\f{\lambda_j(\lambda_k-\lambda_i)+\lambda_i(\lambda_j-\lambda_i)}{\lambda_i^2\lambda_j\lambda_k}\\
    =&\f{F^2}{(n-1)^2}\sum_{j,k}\left(\f{\lambda_k-\lambda_i}{\lambda_i^2\lambda_k}+\f{\lambda_j-\lambda_i}{\lambda_i\lambda_j\lambda_k}\right).
\end{align*}
Hence we obtain 
\[
\abs{F^{ii}-\f1{(n-1)^2}}\leq\f{F^2}{(n-1)^2}|\mathring{A}|\sum_{j,k}\f{\lambda_j+\lambda_k}{\lambda_i^2\lambda_j\lambda_k}\leq C\sigma^{-n-\delta+1}.
\]
On the other hand, 
\begin{equation}\label{e-F-1}
    \abs{F^{kl}_{,l}}=\abs{F^{kl,rs}h_{rs,l}}\leq C\abs{d^2F}\abs{\n A}\leq C\sigma^{-n-\delta}.
\end{equation}
And then combined with Holder inequality, we get 
\begin{align*}
    \int_Mu\mathcal{L}u\leq&-\f1{(n-1)^2}\int_M\abs{\n u}^2+C\sigma^{-n-\delta}\int_M\abs{\n u}\abs{u}+C\int_M\sigma^{-n-\delta+1}\abs{\n u}^2\\
    \leq&-\f1{(n-1)^2}\int_M\abs{\n u}^2+C\sigma^{-n-\delta+1}\int_M\abs{\n u}^2+C\sigma^{-n-\delta-1}\int_Mu^2\\
    =&-\left(\f1{(n-1)^2}-C\sigma^{-n-\delta+1}\right)\int_M\abs{\n u}^2+C\sigma^{-n-\delta-1}\int_Mu^2\\
    \leq&-\left(\f1{(n-1)^2}-C\sigma^{-n-\delta+1}\right)\left(\f{n-1}{\sigma^2}-2\f{(n-1)m}{(n-2)\sigma^n}-C\sigma^{-n-1}\right)\int_Mu^2\\
    \leq&-\left(\f1{(n-1)\sigma^2}-\f {2m}{(n-1)(n-2)\sigma^n}-C\sigma^{-n-1}\right)\int_Mu^2.
\end{align*}
To show the second inequality, we note that $\abs{F^{ij}_{,j}}\leq C\sigma^{-n-\delta}$. We compute
\begin{align*}
    \int_MF^{ij}u_iu_j&=-\int_Mu\mathcal{L}u-\int_MF^{ij}_{,j}u_iu\\
    &\geq\left(\f1{(n-1)^2}-C\sigma^{-n-\delta+1}\right)\int_M\abs{\n u}^2-C\sigma^{-n-\delta-1}\int_Mu^2\\
    &\geq\left(\f1{(n-1)\sigma^2}-\f{2m}{(n-1)(n-2)\sigma^n}-C\sigma^{-n-1}\right)\int_Mu^2.
\end{align*}
\end{proof}
In the next lemma, we will give the estimate of the curvature $F$, which indicates that $F$ is approaching a constant up to polynomial perturbation.
\begin{lem}\label{fF}
    $\abs{F-f}\leq C\sigma^{-n-\delta}$.
\end{lem}
    \begin{proof}
        By the definition, $F_i=F^{kl}h_{kl,i}$. We have
        \[
        \abs{\n F}\leq\abs{dF}\abs{\n A}\leq C\sigma^{-n-\delta-1}.
        \]
        On the other hand, we have 
        \[
        Diam(M)\leq C\int_M\abs{H}^{n-2}\leq C\sigma^{2-n}\sigma^{n-1}\leq C\sigma.
        \]
        It then follows that 
        \[
        \abs{F-f}\leq\max_M\abs{\n F}Diam(M)\leq C\sigma^{-n-\delta}.
        \]
    \end{proof}
    We are now in a position to give the exponential decay of the flow. We recall the evolution equation in \cite{Guilisun} 
    \[
    \f\pa{\pa t}F=\mathcal{L}F-(f-F)(F^{kl}h_{mk}h_{ml}-F^{kl}\bar{R}_{\nu k\nu l}),
    \]
    and from Proposition \ref{prop2.2},
    \[
    \lambda_i=\f1\sigma-\f{(n-1)m}{n-2}\sigma^{1-n}+O(\sigma^{-n}),
    \]
    \[
    \abs{F^{ii}-\f1{(n-1)^2}}\leq C\sigma^{-n-\delta+1}.
    \]
    We then calculate as 
    \begin{align*}
        &F^{kl}h_{mk}h_{ml}\\
        =&\left(\f1{(n-1)^2}+O(\sigma^{-n-\delta+1})\right)(n-1)\left(\f1\sigma-\f{(n-1)m}{n-2}\sigma^{1-n}+O(\sigma^{-n})\right)^2\\
        =&\f1{(n-1)\sigma^2}-\f{2m}{(n-2)\sigma^n}+O(\sigma^{-n-1}),
    \end{align*}
    and by Lemma \ref{lem1.1} and Lemma \ref{lem1.2},
    \begin{align*}
        &F^{kl}\bar{R}_{\nu k\nu l}\\
        =&F^{kl}\Big(\f1{n-2}(-\bar{R}_{\nu\nu}g_{kl}-\bar{R}_{kl})+O(\sigma^{-n-\delta-1})\Big)\\
        =&\f1{n-2}\Big(\f1{(n-1)^2}+O(\sigma^{-n-\delta+1})\Big)\left(\f{(n-1)^2(n-2)m}{\sigma^n}-\f{(n-1)(n-2)m}{\sigma^n}+O(\sigma^{-n-1-\delta})\right)\\
        =&\f{(n-2)m}{(n-1)\sigma^n}+O(\sigma^{-n-1}).
    \end{align*}
    It follows that 
    \begin{equation}\label{Rvv}
    F^{kl}h_{mk}h_{ml}-F^{kl}\bar{R}_{\nu k\nu l}=\f1{(n-1)\sigma^2}-\f{m(n-1)^2+m}{(n-1)(n-2)\sigma^{n}}+O(\sigma^{-n-1}).
    \end{equation}
     By Lemma \ref{mu1} and Lemma \ref{fF} and \eqref{Rvv}, we obtain that
    \begin{align*}
        &\f d{dt}\int_{M_t}(f-F)^2d\mu_t\\
        =&\int_{M_t}\left[-2(f-F)\mathcal{L}F+2(f-F)^2(F^{kl}h_{mk}h_{ml}-F^{kl}\bar{R}_{\nu k\nu l})+H(f-F)^3\right]d\mu_t\\
        \leq&-2\Big(\f1{(n-1)\sigma^2}-\f{2m}{(n-1)(n-2)\sigma^n}-C\sigma^{-n-1}\Big)\int_{M_t}(f-F)^2\\
        &+2\Big(\f1{(n-1)\sigma^2}-\f{m(n-1)^2+m}{(n-1)(n-2)\sigma^n}+C\sigma^{-n-1}\Big)\int_{M_t}(f-F)^2\\
        &+\int_{M_t}H(f-F)^3\\
        =&\Big(\f{-2nm}{(n-1)\sigma^n}+C\sigma^{-n-1}\Big)\int_{M_t}(f-F)^2-\int_{M_t}H(F-f)^3\\
        \leq&\Big(\f{-2nm}{(n-1)\sigma^n}+C\sigma^{-n-1}\Big)\int_{M_t}(f-F)^2+C\sigma^{-n-1-\delta}\int_{M_t}(f-F)^2\\
        \leq&-\f{2m}{\sigma^n}\int_{M_t}(f-F)^2d\mu_t.
    \end{align*}
This implies that
\begin{equation}\label{expo}
    \int_{M_t}(f-F)^2d\mu_t\leq e^{-\frac{2m}{\sigma^n}t}\int_{M_0}(f-F)^2d\mu_0.
\end{equation}
Exponential convergence follows by standard interpolation inequalities,  completing the proof of Theorem \ref{main}.

\vspace{.2in}

\section{Existence of the constant harmonic mean curvature foliation}

\vspace{.1in}

This section addresses the existence of constant harmonic mean curvature hypersurface and consequently a version of "center of gravity" can be defined. We consider the smooth operator $\mathscr{F}: C^3(S^{n-1}, N)\rightarrow C^1(S^{n-1})$ which assigns to each embedding $\phi: S^{n-1}\rightarrow N$ the harmonic mean curvature $\mathscr{F}(\phi)$ of the hypersurface $\Sigma=\phi(S^{n-1})$. 
Given a variation vector field $V$ on a constant harmonic mean curvature hypersurface $\Sigma$, we could compute directly that the first variation of the $\mathscr{F}$ operator at $\Sigma$ in direction $V$ is given by
\[
d\mathscr{F}(\phi)\cdot V=-\left(\mathcal{L}+F^{ij}h_{jk}h_{ik}-F^{ij}\bar{R}_{\nu i\nu j}\right)\<V,\nu\>,
\]
where $\phi:S^{n-1}\rightarrow{N}$ is the embedding of the constant harmonic mean curvature hypersurface.
Define the linearized harmonic mean curvature operator $L$ on $\Sigma$
\begin{equation}\label{e-L}
L=-\left(\mathcal{L}+F^{ij}h_{jk}h_{ik}-F^{ij}\bar{R}_{\nu i\nu j}\right).
\end{equation}
Since the operator $L$ is not self-adjoint, many useful tools in functional analysis and PDE cannot be applied to it. To proceed further, we consider the adjoint operator $L^*$ of $L$. By direct computation, we see that
\begin{equation}\label{e-adjoint-L}
    L^*u=Lu-2F^{ij}_{,i}u_j-F^{ij}_{,ij}u.
\end{equation}
Then the operator
\begin{eqnarray}\label{e-S}
    Su&:=&\frac{1}{2}(L+L^*)u=Lu-F^{ij}_{,i}u_j-\frac{1}{2}F^{ij}_{,ij}u\nonumber\\
    &=& -\left(\mathcal{L}u+F^{ij}h_{jk}h_{ik}u-F^{ij}\bar{R}_{\nu i\nu j}u+F^{ij}_{,i}u_j+\frac{1}{2}F^{ij}_{,ij}u\right)\nonumber\\
    &=& -(F^{ij}u_j)_{,i}-\left(F^{ij}h_{jk}h_{ik}-F^{ij}\bar{R}_{\nu i\nu j}+\frac{1}{2}F^{ij}_{,ij}\right)u
\end{eqnarray}
is a self-adjoint operator.

\iffalse
\begin{defn}
    A smooth hypersurface $\Sigma$ in $N$ is called strictly harmonic stable if the linearized harmonic mean curvature operator $L$ satisfies 
    \[\mu_0:=\inf \left\{\int_{\Sigma}uLud\mu: ||u||_{L^2}=1, \int_{\Sigma}ud\mu=0\right\}>0.\]
\end{defn}
\fi

Denote 
    \[\mu_0:=\inf \left\{\int_{\Sigma}uSud\mu: ||u||_{L^2}=1, \int_{\Sigma}ud\mu=0\right\}>0.\]
We have the following lower bound for $\mu_0$:

\begin{lem}\label{strc}
    Let $\Sigma_{\sigma}$ be the constant harmonic mean curvature hypersurface constructed in Theorem \ref{main}. Then for $\sigma\geq\sigma_0$, we have
    \[ \mu_0\geq \frac{nm}{n-1}\sigma^{-n}-c\sigma^{-n-1}. \]
\end{lem}

\begin{proof}
We denote $\Sigma_{\sigma}$ by $\Sigma$ for any $\sigma\geq\sigma_0$.
By \eqref{Rvv}, we see that 
Using Lemma \ref{lem1.9} and Proposition \ref{3.10}, we also have
\begin{equation}\label{e-F-ij}
F^{ij}_{,ij}=O(\sigma^{-n-1-\delta}).
 \end{equation} 
For any $u$ with $||u||_{L^2}=1, \int_{\Sigma}ud\mu=0$, by Lemma \ref{mu1}, we have
\begin{equation*}
    \begin{split}
        \int_{\Sigma}uSud\mu=& \int_{\Sigma}\left[-u\mathcal{L}u-(F^{ij}h_{jk}h_{ik}-F^{ij}\bar{R}_{\nu i\nu j})u^2-\left(F^{ij}_{,i}u_ju+\frac{1}{2}F^{ij}_{,ij}u^2\right)\right]d\mu\\
        \geq&\frac{nm}{(n-1)\sigma^n}-c\sigma^{-n-1}.
    \end{split}
\end{equation*}
    
\end{proof}

\iffalse
\red{The next lemma is essentially to show that the linearized operator is invertible and the decay estimate. However, the proof employ the eigenfunction decomposition, which is not valid for a non-symmetric operator( here is the case, since our linearized operator is not symmetric, although it can be seen as a perturbation of the Laplacian operator).}  We have to employ another method to give the desired result. (However, yes there is another "but", \blue{the eigenvalues are actually "larger" than the first eigenvalue, in the sense that the real part of the eigenvalue is larger than the first eigenvalue, this fact indeed imply that $L$ is invertible and also a suitable estimate of the norm of the inverse.} 

\fi
A significant consequence of the above lemma is that the operator $S$ is invertible, as the following lemma says,
\begin{lem}\label{inv}
   Let $\Sigma_{\sigma}$ be the constant harmonic mean curvature hypersurface constructed in Theorem \ref{main}. For $\sigma\geq\sigma_0$, the operator $S$ is invertible, and $|S^{-1}|\leq cm^{-1}\sigma^{n}$.
\end{lem}

\begin{proof}
     We denote $\Sigma_{\sigma}$ by $\Sigma$ for any $\sigma\geq\sigma_0$. Let $\eta_0$ be the lowest eigenvalue of $S$ without constraints, i.e.
    \begin{align*}
        \eta_0=\inf_{\{u: ||u||_{L^2}=1\}}\int_{\Sigma}uSud\mu%=\inf_{\{u: ||u||_{L^2}=1\}}\int_{\Sigma}uLud\mu. 
    \end{align*}
    By the proof of Lemma \ref{mu1}, for $u$ with $||u||_{L^2}=1$, there holds, for $\sigma\geq\sigma_0 $
    \begin{align*}
        \int_{\Sigma}u\mathcal{L}ud\mu
        \leq  -\left(\frac{1}{(n-1)^2}-c\sigma^{-n-\delta+1}\right)\int_{\Sigma}|\n u|^2d\mu+c\sigma^{-n-\delta-1}.
    \end{align*}
   By \eqref{Rvv}, we have
   \begin{equation*}
    \begin{split}
        &\int_{\Sigma}uSud\mu\\
        =& \int_{\Sigma}\left[-u\mathcal{L}u-(F^{ij}h_{jk}h_{ik}-F^{ij}\bar{R}_{\nu i\nu j})u^2-\left(F^{ij}_{,i}u_ju+\frac{1}{2}F^{ij}_{,ij}u^2\right)\right]d\mu\\
        \geq&\left(\frac{1}{(n-1)^2}-c\sigma^{-n-\delta+1}\right)\int_{\Sigma}|\n u|^2d\mu-\f1{(n-1)\sigma^2}+\f{(n^2-2n+2)m}{(n-1)(n-2)\sigma^n}+O(\sigma^{-n-1}).
    \end{split}
\end{equation*}
Hence,
   \[ \eta_0\geq-\frac{1}{(n-1)\sigma^2}+\frac{(n^2-2n+2)m}{(n-1)(n-2)\sigma^n}-c\sigma^{-n-1}.  \]
    On the other hand, if we replace $u$ by a constant, we obtain the reverse inequality. Therefore, 
    \begin{equation*}
        \eta_0=-\frac{1}{(n-1)\sigma^2}+\frac{(n^2-2n+2)m}{(n-1)(n-2)\sigma^n}+O(\sigma^{-n-1}).
    \end{equation*}
    Let $h_0$ be the corresponding eigenfunction of $\eta_0$
    \[ Sh_0=\eta_0h_0. \]
    Let $\bar{h}_0=\fint_{\Sigma}h_0d\mu$ be the mean value of $h_0$.
    Multiplying the above identity with $(h_0-\bar{h}_0)$ and integrating it over $\Sigma$ to obtain
    \begin{align*}
        & -\int_{\Sigma}(h_0-\bar{h}_0)\mathcal{L}(h_0-\bar{h}_0)d\mu\\
        =& \int_{\Sigma}\left(\eta_0+F^{ij}h_{jk}h_{ik}-F^{ij}\bar{R}_{\nu i\nu j}+\frac{1}{2}F^{ij}_{,ij}\right)(h_0-\bar{h}_0)^2d\mu\\
        &+\int_{\Sigma}\left(\eta_0+F^{ij}h_{jk}h_{ik}-F^{ij}\bar{R}_{\nu i\nu j}+\frac{1}{2}F^{ij}_{,ij}\right)\bar{h}_0(h_0-\bar{h}_0)d\mu\\
        &+\int_{\Sigma}F^{ij}_{,i}h_{0,j}(h_0-\bar{h}_0)d\mu\\
        =& \int_{\Sigma}\left(\eta_0+F^{ij}h_{jk}h_{ik}-F^{ij}\bar{R}_{\nu i\nu j}\right)(h_0-\bar{h}_0)^2d\mu\\
        &+\int_{\Sigma}\left(\eta_0+F^{ij}h_{jk}h_{ik}-F^{ij}\bar{R}_{\nu i\nu j}+\frac{1}{2}F^{ij}_{,ij}\right)\bar{h}_0(h_0-\bar{h}_0)d\mu.
    \end{align*}
    From Lemma \ref{mu1}, the left hand side is bounded below by
    \[ -\int_{\Sigma}(h_0-\bar{h}_0)\mathcal{L}(h_0-\bar{h}_0)d\mu\geq \left(\f1{(n-1)\sigma^2}-\f{2m}{(n-1)(n-2)\sigma^n}-c\sigma^{-n-1}\right)\int_M(h_0-\bar{h}_0)^2d\mu. \]
    Together with 
    \[ \eta_0+F^{ij}h_{jk}h_{ik}-F^{ij}\bar{R}_{\nu i\nu j}=O(\sigma^{-n-1}) \]
    and (\ref{e-F-ij}), we obtain, for $\sigma\geq\sigma_0$,  
    \begin{equation}\label{h0}
        ||h_0-\bar{h}_0||_{L^2}\leq c\sigma^{1-n}|\bar{h}_0||\Sigma|^{\frac{1}{2}}.
    \end{equation}
    In particular, $\bar{h}_0\not=0$. Let $\eta_1$ be the next eigenvalue of $S$ with corresponding eigenfunction $h_1$.
    \begin{comment}
        here in general, $\eta_1$ is not a real number. But we have the following lemma.
    \begin{lem}
    Principal eigenvalue for non-symmetric elliptic operators,
    \begin{itemize}
        \item i)There exists a real eigenvalue $\lambda_1$ for the operator L,  taken with zero mean value integral, such that if $\lambda\in\mathbb{C}$ is any other eigenvalue,  we have
$$\mathrm{Re}(\lambda)\geq\lambda_1.$$
\item ii) There exists a corresponding eigenfunction $w_1$, which is positive within U.
\item iii) The eigenvalue $\lambda _1$ is simple; that is,  if u is any solution,  then u is a multiple of $w_1$.
  \end{itemize}
    \end{lem}
\end{comment}
    Let $\bar{h}_1=|\Sigma|^{-1}\int_{\Sigma}h_1d\mu$ be the mean value of $h_1$. 
    Note that
    \[ 0=\int_{\Sigma}h_0h_1d\mu=\int_{\Sigma}(h_0-\bar{h}_0)(h_1-\bar{h}_1)d\mu+\int_{\Sigma}\bar{h}_0h_1d\mu. \]
    By H\"older inequality, we get
    \[ \left|\int_{\Sigma}h_1d\mu\right|\leq |\bar{h}_0|^{-1}||h_0-\bar{h}_0||_{L^2}||h_1-\bar{h}_1||_{L^2}. \]
    Then, by \eqref{h0}, there holds
    \begin{equation}\label{h1}
        |\bar{h}_1|\leq c\sigma^{1-n}|\Sigma|^{-\frac{1}{2}}||h_1-\bar{h}_1||_{L^2}.
    \end{equation}
Multiplying $Sh_1=\eta_1 h_1$ with $(h_1-\bar{h}_1)$ and integrating it over $\Sigma$,  we obtain by Lemma \ref{strc} and \eqref{h1},
\begin{align*}
&\left(\frac{nm}{(n-1)\sigma^{n}}-c\sigma^{-n-1}\right)\int_{\Sigma}(h_1-\bar{h}_1)^2d\mu\\
    \leq& \int_{\Sigma}(h_1-\bar{h}_1)S(h_1-\bar{h}_1)d\mu\\
    =& \eta_1\int_{\Sigma}(h_1-\bar{h}_1)^2d\mu+\bar{h}_1\int_{\Sigma}(h_1-\bar{h}_1)\left(F^{ij}h_{jk}h_{ik}-F^{ij}\bar{R}_{\nu i\nu j}+\frac{1}{2}F^{ij}_{,ij}\right)d\mu\\
    \leq &\eta_1\int_{\Sigma}(h_1-\bar{h}_1)^2d\mu
    +c\sigma^{-n-1}|\bar{h}_1|\int_{\Sigma}|h_1-\bar{h}_1|d\mu\\
    \leq &\eta_1\int_{\Sigma}(h_1-\bar{h}_1)^2d\mu+c\sigma^{-2n}\int_{\Sigma}(h_1-\bar{h}_1)^2d\mu.
\end{align*}
Therefore, for $\sigma\geq\sigma_0$, 
\[ \eta_1\geq cm\sigma^{-n}. \]
Hence, we show that $S$ is invertible, and $|S^{-1}|\leq cm^{-1}\sigma^n$.
\end{proof}

We are now ready to use Lemma 4.1 in \cite{Guilisun} to show that $\{\Sigma_{\sigma}\}$ form a foliation. 

\begin{thm}\label{local}
    There is $\sigma_0$ depending only on $c_0$ such that for all $\sigma\geq\sigma_0$, the constant harmonic mean curvature hypersurfaces $\Sigma_{\sigma}$ constructed in Theorem \ref{main} constitute a proper foliation of $N\backslash B_{\sigma_0}(0)$. 
\end{thm}

\begin{proof}
    Let $\Sigma_{\sigma}$ be the family of constant harmonic mean curvature hypersurfaces constructed in Theorem \ref{main} with $\sigma\geq\sigma_0$, and $\phi^{\sigma}: S^{n-1}\rightarrow N$ be the embedding for each $\sigma$.

    For any $\sigma_2>\sigma_1\geq\sigma_0$ such that $\Sigma_{\sigma_i}\in \mathcal{B}_{\sigma}(B_1, B_2, B_3), i=1,2$, $\Sigma_{\sigma_2}$ can be represented by a normal variation over $\Sigma_{\sigma_1}$ of the form
    \[ \phi^{\sigma_2}=\phi^{\sigma_1}+u(p)\nu(p), \quad p\in S^{n-1}.
    \]
     We will show that $u$ has a sign; in particular, $u$ cannot be zero.
    In the following, we denote $\Sigma_{\sigma}$ by $\Sigma$.
    
Let $V=\phi^{\sigma_2}-\phi^{\sigma_1}$.
    By the Taylor theorem, there is some $\sigma\in(\sigma_1, \sigma_2)$ such that
    \begin{eqnarray}\label{e-Taylor}
        \mathscr{F}(\phi^{\sigma_2})&=&\mathscr{F}(\phi^{\sigma_1})+Lu+\frac{1}{2}d^2\mathscr{F}(\phi^{\sigma})(V, V)\nonumber\\
    &=&\mathscr{F}(\phi^{\sigma_1})+Su+(L-S)u+\frac{1}{2}d^2\mathscr{F}(\phi^{\sigma})(V, V),
    \end{eqnarray}
    where $\mathscr{F}(\phi^{\sigma_2})$ and $\mathscr{F}(\phi^{\sigma_1})$ are constants.
    Note that the second variation of the harmonic mean curvature operator involves the second fundamental form and higher derivatives of the metric, this leads to
    \begin{equation}\label{ddf}
        ||d^2\mathscr{F}(V,V)||\leq c\sigma^{-3}||V||^2\leq c\sigma^{-3}||u||_{C^{2}}.
    \end{equation}  
    We also have from the definition of $S$ that
       \begin{equation}\label{e-S-L}
        |(L-S)u|=\frac{|(L-L^*)u|}{2}=\left|F^{ij}_{,i}u_j+\frac{1}{2}F^{ij}_{,ij}u\right|\leq c\sigma^{-n-\delta}||u||_{C^{1}}.
    \end{equation}  
Hence, $u$ satisfies the following elliptic equation
\begin{equation}\label{lu}
    Su=\mathscr{F}(\phi^{\sigma_2})-\mathscr{F}(\phi^{\sigma_1})+E_1,
\end{equation}  
with the error term satisfying $|E_1|\leq c\sigma^{-3}||u||_{C^{2}}$.
We will employ Huang's approach of Theorem 3.9 in \cite{Huang}.
By integrating \eqref{lu} over $S^{n-1}$, 
\[ \mathscr{F}(\phi^{\sigma_2})-\mathscr{F}(\phi^{\sigma_1})=
        \frac{1}{|S^{n-1}|}\int_{S^{n-1}}Sud\mu-\frac{1}{|S^{n-1}|}\int_{S^{n-1}}E_1d\mu. \]
By the definition of $S$, we get
\begin{equation*}
    \begin{split}
        &\frac{1}{|S^{n-1}|}\int_{S^{n-1}}Sud\mu\\
        =& -\frac{1}{|S^{n-1}|}\int_{S^{n-1}}\left(\mathcal{L}u+F^{ij}h_{ik}h_{jk}u-F^{ij}\bar{R}_{\nu i\nu j}u+F^{ij}_{,i}u_j+\frac{1}{2}F^{ij}_{,ij}u\right)d\mu\\
        =& -\frac{1}{|S^{n-1}|}\int_{S^{n-1}}\left(\mathcal{L}u+F^{ij}h_{ik}h_{jk}u-F^{ij}\bar{R}_{\nu i\nu j}u-\frac{1}{2}F^{ij}_{,ij}u\right)d\mu. 
    \end{split}
\end{equation*}
Note that 
\[ \mathcal{L}u=F^{ij}u_{ij}=\frac{\Delta u}{(n-1)^2}+O(|\mathring{A}||A|^{-1})u_{ij}. \]
So
\begin{equation*}
    \begin{split}
        &\frac{1}{|S^{n-1}|}\int_{S^{n-1}}Sud\mu\\
        \leq& \left(-\f1{(n-1)\sigma^2}+\f{(n^2-2n+2)m}{(n-1)(n-2)\sigma^n}+O(\sigma^{-n-1})\right)\bar{u}+c\sigma^{-n-\delta}\norm{u}_{C^1},
    \end{split}
\end{equation*}
where $\bar{u}=\frac{1}{|S^2|}\int_{S^2}ud\mu$ is the mean value of $u$ and the last inequality is obtained by integration by parts.
It follows,
\begin{equation}\label{fphi}
    \mathscr{F}(\phi^{\sigma_2})-\mathscr{F}(\phi^{\sigma_1})=-\frac{1}{(n-1)\sigma^2}\bar{u}+E_2,
\end{equation}
where 
\[ |E_2|\leq c\sigma^{-3}||u||_{C^2}. \]
Now we decompose $u=h_0+u_0$, where $h_0$ is the lowest eigenfunction of $S$ and $\int_{S^2}h_0u_0d\mu=0$.

Due to De Giorgi-Nash-Moser iteration and \eqref{h0}, there holds the estimate
\[ \sup|h_0-\bar{h}_0|\leq c\sigma^{1-n}|\bar{h}_0|. \]
We note that $(h_0-\bar{h}_0)$ satisfies the following equation
\[ S(h_0-\bar{h}_0)=\eta_0(h_0-\bar{h}_0)+\left(\eta_0+F^{ij}h_{ik}h_{jk}-F^{ij}\bar{R}_{\nu i\nu j}+\frac{1}{2}F^{ij}_{,ij}\right)\bar{h}_0. \]
And hence we obtain 
\begin{equation}\label{hal}
\begin{split}
    &||h_0-\bar{h}_0||_{C^{0, \alpha}}\\
    \leq& c\left(\sigma^{-\alpha}\sup_{S^{n-1}}|h_0-\bar{h}_0|+|\bar{h}_0|\norm{\eta_0+F^{ij}h_{ik}h_{jk}-F^{ij}\bar{R}_{\nu i\nu j}+\frac{1}{2}F^{ij}_{,ij}}_{L^{n-1}}\right)\\
    \leq & c\sigma^{1-n-\alpha}|\bar{h}_0|,
\end{split}     
\end{equation}
for some $\alpha\in(0, 1)$.

Since $u_0=u-h_0$, \eqref{lu} and \eqref{fphi} imply
\begin{equation*}
    \begin{split}
        Su_0=& Su-Sh_0=-\frac{1}{(n-1)\sigma^2}\bar{u}+E_2+E_1-\eta_0h_0\\
        =& \frac{1}{(n-1)\sigma^2}(h_0-\bar{h}_0)-\frac{1}{(n-1)\sigma^2}\bar{u}_0+E_3+E_2+E_1,
    \end{split}
\end{equation*}
where $|E_3|\leq c\sigma^{-n}|\bar{h}_0|$.
Because $S$ has no kernel,
the Schauder estimate and \eqref{hal} gives
\begin{equation}\label{u0c2}
    \begin{split}
        \norm{u_0}_{C^{2, \alpha}}\leq & c(\sigma^{-2}\norm{h_0-\bar{h}_0}_{C^{0, \alpha}}+\sigma^{-2}|\bar{u}_0|+\sigma^{-3}\norm{u}_{C^{2, \alpha}}+\sigma^{-n}|\bar{h}_0|)\\
        \leq & c(\sigma^{-1-n-\alpha}|\bar{h}_0|+\sigma^{-2}|\bar{u}_0|+\sigma^{-3}(\norm{u_0}_{C^{2, \alpha}}+\norm{h_0}_{C^{2, \alpha}})+\sigma^{-n}|\bar{h}_0|).
    \end{split}
\end{equation}
Again since $h_0$ satisfies $Sh_0=\eta_0h_0$ and $\eta_0=O(\sigma^{-2})$, inequality \eqref{hal} leads to
\[ \norm{h_0}_{C^{2, \alpha}}\leq c\sigma^{-2}\norm{h_0}_{C^{0, \alpha}}\leq c\sigma^{-2}(\norm{h_0-\bar{h}_0}_{C^{0, \alpha}}+|\bar{h}_0|)\leq c\sigma^{-2}|\bar{h}_0|. \]
Inserting this into \eqref{u0c2}, we arrive at
\begin{equation}\label{u0}
    \norm{u_0}_{C^{2, \alpha}}\leq c(\sigma^{-\min\{5,n\}}|\bar{h}_0|+\sigma^{-2}|\bar{u}_0|),
\end{equation}
for $\sigma\geq\sigma_0$. 
Notice the fact that 
\[0=\int_{S^{n-1}}h_0u_0d\mu=\int_{S^{n-1}}u_0(h_0-\bar{h}_0)+\bar{h}_0\int_{S^{n-1}}u_0,\] 
this implies 
\[ |\bar{u}_0|\leq |\bar{h}_0|^{-1}\sup |h_0-\bar{h}_0|\sup |u_0|\leq c\sigma^{1-n}|u_0|_{C^{2, \alpha}}. \]
Then combining with \eqref{u0} we arrive at the following inequality
\[ \norm{u_0}_{C^{2, \alpha}}\leq c\sigma^{-\min\{5,n\}}|\bar{h}_0|. \]
Therefore, we obtain
\[ |u-\bar{h}_0|\leq |u_0|+|h_0-\bar{h}_0|\leq c\sigma^{-\min\{5,n-1\}}|\bar{h}_0|,\]
which shows that $u$ has a sign since $\bar{h}_0$ has a sign.
Finally, an application of Lemma 4.1 in \cite{Guilisun} completes the proof.
\end{proof}

We will now define the ``center of gravity" as the constant harmonic mean curvature foliation approaches infinity. Following the way of Corvino-Wu in \cite{Corvino}, we will show that the center of mass defined by the constant harmonic mean curvature foliation is equivalent to the ADM center of mass.
We first present some definitions and necessary propositions.

The second derivative of the second fundamental form can be similarly estimated as Proposition \ref{3.10}, see the appendix in \cite{Corvino}. Together with Lemma \ref{lem1.4}, we can derive the following lemma.
\begin{lem}\label{2of}
    Let $\Sigma_{\sigma}$ be the constant harmonic mean curvature hypersurface constructed in Theorem \ref{main} with $\delta\geq 0$. Then, for all $t\geq0$ and $\sigma\geq\sigma_0$ large enough,
    \[ |\n^2F|\leq c\sigma^{-n-2-\delta}, \]
    where $c$ depends on $c_0, B_1, B_2, B_3$.
\end{lem}

\begin{defn}
    Assume $(M,g)$ be an asymptotically Schwarzschild manifold with metric $g_{ij}=\left(1+\frac{m}{2|x|^{n-2}} \right)^{\frac{4}{n-2}}\delta_{ij}+O_2(|x|^{-2})$, the ADM center of mass is defined by 
{\footnotesize   \[ C^k=\frac{1}{2m(n-1)\omega_{n-1}}\lim_{R\to\infty}\int_{|x| = R}\left[\sum_{i}x^{k}\left(g_{ij,i}-g_{ii,j}\right)v_{e}^{j}-\sum_{i}\left(g_{ik}v_{e}^{i}-g_{ii}v_{e}^{k}\right)\right]d\mu_{e},\]}
    where $\nu_e$ is the normal vector with respect to the Euclidean metric $\delta$, $k=1, \cdots, n$.
\end{defn}

Consider an exterior region in an asymptotically flat manifold, and assume that there is an asymptotically flat chart in which $g_{ij}=\left(1+\frac{m}{2|x|^{n-2}}+\frac{B_{k}x^{k}}{|x|^{n}}\right)^{\frac{4}{n-2}}\delta_{ij}+O_{5}\left(|x|^{-n}\right)$ with $m>0$. By direct calculation, the ADM center of mass is 
\begin{equation}\label{center}
    C^k=\frac{2B_k}{(n-2)m}. 
\end{equation}

\begin{defn}\label{center}
    Let $\Sigma_{\sigma}$ be the family of surfaces constructed in Theorem \ref{local} and $\phi^{\sigma}$ be the position vector.
     The center of gravity of $\Sigma_{\sigma}$ is defined as
    \[ C_{HM}=\lim_{\sigma\rightarrow\infty}\frac{\int_{\Sigma_{\sigma}}\phi^{\sigma} d\mu_e}{\int_{\Sigma_{\sigma}}d\mu_e}. \]
\end{defn}

We note that translating the coordinates by a vector $a^k$ shifts both the ADM center of mass and the center $C_{HM}$ by $a^k$. Thus, we can adjust the coordinates to set the ADM center of mass vector to be zero, which results in $B_k=0$ in the chart by \eqref{center}.   

\begin{thm}
    Let $(N,\bar{g})$ be an asymptotically flat manifold whose exterior region admits a foliation by constant harmonic mean curvature hypersurfaces, as established in Theorem \ref{local}. Suppose there exists an asymptotically flat coordinate chart in which the metric satisfies 
    \[\bar{g}_{ij}=\left(1+\frac{m}{2|x|^{n-2}} \right)^{\frac{4}{n-2}}\delta_{ij}+O_{5}\left(|x|^{-n}\right)\]
    with $m>0$.
    Then $C_{HM}=0$.
\end{thm}

\begin{proof}
    By the Sobolev embedding inequality and interpolation inequality, we have for $n-1<q<r$ and $\lambda\in (0,1)$ that 
\begin{align}\label{inter}
    ||F-f||_{C^0}&\leq c\sigma^{-\frac{n-1}{q}}(\sigma||\n F||_{L^q}+||F-f||_{L^q})\nonumber\\
    &\leq c\sigma^{-\frac{n-1}{q}}(\sigma||\n F||_{L^q}+||F-f||_{L^r}^{1-\lambda}||F-f||_{L^2}^{\lambda}), 
\end{align} 
where $\frac{1}{q}=\frac{1-\lambda}{r}+\frac{\lambda}{2}$ and norms are calculated on $\Sigma_t$ which is a solution to \eqref{flow} contained in $\mathcal{B}_{\sigma}(B_1, B_2, B_3)$.
At $t=0$, we are working in a coordinate system that $g_{ij}-g^S_{ij}=O_5(|x|^{-n})$, with $g^S_{ij}$ being the Schwarzschild metric. 
Thus, the difference between the second fundamental form of $\Sigma_t$ is 
\[ h_{ij}-h^S_{ij}=O(|x|^{-n-1}), \]
    which implies that $|F-f|=O(\sigma^{-n-1})$.
From \eqref{expo}, for all $t\geq0$,
\[ \int_{\Sigma_t}(F-f)^2d\mu_t\leq C\sigma^{-n-3}e^{-\frac{2mt}{\sigma^n}}. \]
To derive a pointwise bound from \eqref{inter}, we utilize the interpolation inequality from \cite{Hamilton} 
to estimate, for $\frac{1}{p}+\frac{1}{2}=\frac{2}{q}$,
\begin{equation*}
    \begin{split}
        \sigma^{1-\frac{n-1}{q}}||\n F||_{L^q}
       \leq & (q+n-3)^{\frac{1}{2}}\sigma^{1-\frac{n-1}{q}}||\n^2F||_{L^p}^{\frac{1}{2}}||F-f||_{L^2}^{\frac{1}{2}} \\
       \leq & C \sigma^{1-\frac{n-1}{q}}(\sigma^{-n-2}\sigma^{\frac{n-1}{p}})^{\frac{1}{2}}\cdot(\sigma^{-n-3}e^{-\frac{2mt}{\sigma^n}})^{\frac{1}{4}}\\
       = & C\sigma^{-n-\frac{1}{2}}e^{-\frac{mt}{2\sigma^n}},
    \end{split}
\end{equation*} 
and 
\begin{equation*}
    \begin{split}
        &\sigma^{-\frac{n-1}{q}}||F-f||_{L^r}^{1-\lambda}||F-f||_{L^2}^{\lambda}\\
        \leq & C\sigma^{-\frac{n-1}{q}}\cdot\left(\sigma^{-n-1}\sigma^{\frac{n-1}{r}}\right)^{1-\lambda}\cdot\left(\sigma^{-n-3}e^{-2mt\sigma^{-n}}\right)^{\frac{\lambda}{2}}\\
        = & C\sigma^{-n-1}e^{-\frac{mt}{2\sigma^n}}.
    \end{split}
\end{equation*}
In the last inequality, we have selected that $\lambda=\frac{1}{2}$.
Therefore, combining with \eqref{inter} we get on $\Sigma_t$ that, for all $t\geq0$ and $\sigma\geq\sigma_0$,
\[ |F-f|\leq C\sigma^{-n-\frac{1}{2}}e^{-\frac{mt}{2\sigma^n}}. \]
By integrating the flow equation $\frac{\partial}{\partial t}\phi^{\sigma}=(f-F)\nu$, and $\partial_t(d\mu_g)=H(f - F)d\mu_g$, we obtain 
\begin{align*}
    |\phi^{\sigma}_{\infty}-\phi^{\sigma}_0|\leq & C\sigma^{-\frac{1}{2}},\\
    |(\phi^{\sigma}_{\infty})^*(d\mu_g)-(\phi^{\sigma}_0)^*(d\mu_g)|\leq & C\sigma^{-\frac{3}{2}}|(\phi^{\sigma}_0)^*(d\mu_g)|.
\end{align*}
  Finally, we have the center of mass 
  \[ C_{HM}=\lim_{\sigma\rightarrow\infty}\frac{\int_{\Sigma_{\sigma}}\phi^{\sigma} d\mu_e}{\int_{\Sigma_{\sigma}}d\mu_e}
  =\lim_{\sigma\rightarrow\infty}\frac{\int_{S^n}\phi^{\sigma}_{\infty}d\mu^{\sigma}}{\int_{S^n}d\mu^{\sigma}}=\lim_{\sigma\rightarrow\infty}O(\sigma^{-\frac{1}{2}})=0,\]
where $d\mu^{\sigma}=(\phi^{\sigma}_{\infty})^*(d\mu_e)$.
\end{proof}

\vspace{.2in}

\section{Uniqueness of the CHMC foliation of dimension 3}

\vspace{.1in}

In this section, we aim to derive the uniqueness of the constant harmonic mean curvature foliation for three-dimensional asymptotically Schwarzschild manifolds under the enhanced decay condition $\delta>0$. Adapting the methodology developed by Huisken and Yau \cite{Huisken1996DefinitionOC}, we prove the following main result.

\begin{thm}\label{thm-uniqueness}
Let $(N^3,\bar g)$ be an asymptotically Schwarzschild space given by (\ref{metr}) with $\delta>0$. Then given $B_1, B_2, B_3>0$, we can choose $\sigma\geq\sigma_0$ such that $\Sigma_{\sigma}$ constructed in Theorem \ref{main} is the only surface with constant harmonic mean curvature $F_{\sigma}$ contained in $\mathcal{B}_{\sigma} (B_1, B_2, B_3)$.
\end{thm}

\begin{proof} We first proceed the computation with general $n$ and restrict ourselves to the case $n=3$ in the end.

Consider the smooth operator 
\begin{align*}
    \mathscr{F}: C^3(S^{n-1}, N^{n})&\to C^1(S^{n-1})\\
    \phi&\mapsto \mathscr{F}(\phi),
\end{align*}
where $\mathscr{F}(\phi)$ is the harmonic mean curvature of the hypersurface given by $\phi$. Recall that the linearization operator of $\mathscr{F}$ is the operator $L$ given by (\ref{e-L}). Now (\ref{e-S}), (\ref{e-F-1}) and (\ref{e-F-ij}), together with Lemma \ref{inv} imply that $L$ is invertible with $|L^{-1}|\leq cm^{-1}\sigma^n$. Here, we used the assumption that $\delta>0$.

\vspace{.1in}

Now suppose there is another hypersurface $\Sigma_{\sigma'}$ with constant harmonic mean curvature $F_{\sigma}$. Let $\phi_0$ and $\phi_1$ denote the embeddings from $S^{n-1}$ to $N^n$, which represent the hypersurfaces $\Sigma_{\sigma}$ and $\Sigma_{\sigma'}$ respectively. Then from Proposition \ref{prop2.1}, we see that there is a vector $\vec{a}$ such that the second surface can be represented by a normal variation of the form
\begin{equation*}
    \phi_1(p)=\phi_0(p)+u(p)\nu(p), \ \ p\in S^{n-1},
\end{equation*}
where $u=\langle \vec{a},\nu\rangle+G$ and $G$ satisfies
\begin{equation*}
\sigma^{n+\delta-2}|G|+\sigma^{n+\delta-1}|\nabla G|+\sigma^{n+\delta}|\nabla^2 G|+\sigma^{n+\delta+1}|\nabla^3 G|\leq c.
\end{equation*}
If we consider the one-parameter family of surfaces $\Sigma_t$ given by
\begin{equation*}
\phi_t=\phi_0+tu\nu,
\end{equation*}
we see that all interpolated hypersurfaces are still in $\mathcal{B}_{\sigma} (B_1, B_2, B_3)$ such that the operator $L$ satisfies the eigenvalue estimates established above on all these hypersurfaces. Since $\mathscr{F}(\phi_0)=\mathscr{F}(\phi_1)=F_{\sigma}$, we get from Taylor's expansion (see (\ref{e-Taylor})) that the
variation vector field $V=\phi_1-\phi_0$ satisfies both
    \begin{equation}\label{e-dF}
        ||d\mathscr{F}(\phi_0)(V)||\leq \sup_{t\in [0,1]}||d^2\mathscr{F}(\phi_t)(V,V)||,
    \end{equation} 
and   
    \begin{equation}\label{e-dF2}
        ||d\mathscr{F}(\phi_0)(V)||\leq \sup_{t\in [0,1]}||(d\mathscr{F}(\phi_t)-d\mathscr{F}(\phi_0)(V)||.
    \end{equation} 
By (\ref{ddf}), we have $||d^2\mathscr{F}(V,V)||\leq c\sigma^{-3}||V||^2$ so that combining with the upper bound for $L^{-1}$
\begin{align*}
    ||V||&=||L^{-1}(LV)||\leq ||L^{-1}||\cdot ||LV||\\
    &\leq  ||L^{-1}||\cdot \sup_{t\in [0,1]}||d^2\mathscr{F}(\phi_t)(V,V)||\\
    &\leq cm^{-1}\sigma^n\cdot c\sigma^{-3}||V||^2
    \leq cm^{-1}\sigma^{n-3}||V||^2
\end{align*}
so that 
\begin{align*}
    ||V||\left(1-cm^{-1}\sigma^{n-3}||V||\right)\leq 0.
\end{align*}
Therefore, if
\begin{align}\label{e-V}
    ||V||<c'm\sigma^{3-n},
\end{align}
for some constant $c'>0$, then we must have $V\equiv0$. By the definition of $u$, we see that if $|\vec{a}|\leq c'm\sigma^{3-n}$ for some constant $c$, then $V\equiv0$.

\vspace{.1in}

To derive the upper bound for $|\vec{a}|$, we note that for all $t\in [0,1]$, 
\begin{align*}
(d\mathscr{F}(\phi_t)(V)=L_t\langle V,\nu\rangle+V^T(F_t),
\end{align*}
where $F_t$ is the harmonic mean curvature of the hypersurface defined by $\phi_t$. Notice that for a hypersurface in $\mathcal{B}_{\sigma} (B_1, B_2, B_3)$, we have
\begin{align*}
||V^T(F_t)||\leq c\sigma^{-n-1-\delta},
\end{align*}
and hence
\begin{align*}
||(d\mathscr{F}(\phi_t)-d\mathscr{F}(\phi_0)(V)||\leq c\sigma^{-n-1}.
\end{align*}
Finally, using the fact that $\langle \vec{a}, \nu\rangle$ is an approximate eigenfunction of $L$ with eigenvalue close to $\frac{nm}{(n-1)\sigma^n}$, we see using (\ref{e-dF2}) that
\begin{align*}
    c\sigma^{-n-1}&\geq\sup_{t\in [0,1]}||(d\mathscr{F}(\phi_t)-d\mathscr{F}(\phi_0)(V)||\geq  Lu\geq c(n)m\sigma^{-n}|\vec{a}|.
\end{align*}
Thus we have
\begin{align*}
  |\vec{a}|\leq cm\sigma^{-1}.
\end{align*}
Now suppose $n=3$. Then if $\sigma$ is sufficiently large, we would have $|\vec{a}|<c'm$ and hence $V\equiv0$. This finishes the proof of the uniqueness. 
\end{proof}

\textbf{Conflict of Interest} The authors have no conflict of interest to declare.

\vspace{.1in}
\textbf{Data availability} The authors declare no datasets were generated or analysed during the current study.

\bibliography{refofmix}

\begin{thebibliography}{Ham82}

\bibitem[CW08]{Corvino}
Justin Corvino and Haotian Wu.
\newblock On the center of mass of isolated systems.
\newblock {\em Classical Quantum Gravity}, 25(8):085008, 18, 2008.

\bibitem[EM13]{Eich-Met}
Michael Eichmair and Jan Metzger.
\newblock Unique isoperimetric foliations of asymptotically flat manifolds in all dimensions.
\newblock {\em Invent. Math.}, 194(3):591--630, 2013.

\bibitem[GLS24]{Guilisun}
Yaoting Gui, Yuqiao Li, and Jun Sun.
\newblock Foliation of constant harmonic mean curvature surfaces in asymptotic schwarzschild spaces.
\newblock {\em arXiv:2412.17024}, 2024.

\bibitem[Ham82]{Hamilton}
Richard~S. Hamilton.
\newblock Three-manifolds with positive {R}icci curvature.
\newblock {\em J. Differential Geometry}, 17(2):255--306, 1982.

\bibitem[Hua09]{Huang09}
Lan-Hsuan Huang.
\newblock On the center of mass of isolated systems with general asymptotics.
\newblock {\em Classical Quantum Gravity}, 26(1):015012, 25, 2009.

\bibitem[Hua10]{Huang}
Lan-Hsuan Huang.
\newblock Foliations by stable spheres with constant mean curvature for isolated systems with general asymptotics.
\newblock {\em Comm. Math. Phys.}, 300(2):331--373, 2010.

\bibitem[Hui86]{Huisken1986ContractingCH}
Gerhard Huisken.
\newblock Contracting convex hypersurfaces in riemannian manifolds by their mean curvature.
\newblock {\em Inventiones mathematicae}, 84:463--480, 1986.

\bibitem[Hui87]{Huisken1987TheVP}
Gerhard Huisken.
\newblock The volume preserving mean curvature flow.
\newblock {\em Journal f{\"u}r die reine und angewandte Mathematik}, 1987:35--48, 1987.

\bibitem[HY96]{Huisken1996DefinitionOC}
Gerhard Huisken and Shing-Tung Yau.
\newblock Definition of center of mass for isolated physical systems and unique foliations by stable spheres with constant mean curvature.
\newblock {\em Inventiones mathematicae}, 124:281--311, 1996.

\bibitem[QT07]{QingTian}
Jie Qing and Gang Tian.
\newblock On the uniqueness of the foliation of spheres of constant mean curvature in asymptotically flat 3-manifolds.
\newblock {\em J. Amer. Math. Soc.}, 20(4):1091--1110, 2007.

\end{thebibliography}
\bibliographystyle{alpha}
    
\end{document}